\documentclass[12pt]{article}
  \usepackage[cp1251]{inputenc}
  \usepackage[english]{babel}
  \usepackage{amsfonts,amssymb,amsmath,amstext,amsbsy,amsopn,amscd,graphicx}
  \usepackage{graphics}
  \usepackage{amsthm}

  \def\thtext#1{
  \catcode`@=11
  \gdef\@thmcountersep{. #1}
  \catcode`@=12}

  \newtheorem{theorem}{Theorem}[section]
  \newtheorem{prop}{Proposition}[section]
  \newtheorem{lem}{Lemma}[section]

 \newcounter{il}[section]
 \renewcommand{\thtext}
 {\thechapter.\arabic{il}}

 \newcounter{il2}[section]
 \renewcommand{\thtext}
 {\thechapter.\arabic{il2}}

 \newenvironment{dfn}{\trivlist \item[\hskip\labelsep{\bf Definition}]
 \refstepcounter{il}{\bf \arabic{section}.\arabic{il}.}}%
 {\endtrivlist}

 \newenvironment{rk}{\trivlist \item[\hskip\labelsep{\bf Remark}]
 \refstepcounter{il2}{\bf \arabic{section}.\arabic{il2}.}}%
 {\endtrivlist}

 \newenvironment{examp}{\trivlist \item[\hskip\labelsep{\bf Example.}]}%
 {\endtrivlist}

 \newenvironment{sproof}{\trivlist \item[\hskip \labelsep{\it Sketch of the proof.}]}%
 {\endtrivlist}

 \catcode`@=11
 \def\.{.\spacefactor\@m}
\def\relaxnext@{\let\next\relax}
\def\nolimits@{\relaxnext@
 \DN@{\ifx\next\limits\DN@\limits{\nolimits}\else
  \let\next@\nolimits\fi\next@}%
 \FN@\next@}
\def\newmcodes@{\mathcode`\'39\mathcode`\*42\mathcode`\."613A%
 \mathcode`\-45\mathcode`\/47\mathcode`\:"603A\relax}
\def\operatorname#1{\mathop{\newmcodes@\kern\z@
 \operator@font#1}\nolimits@}
 \catcode`@=12
 \def\rom#1{{\em#1}}

 \def\({\rom(}
 \def\){\rom)}
 \def\:{\colon}
 \def\corank{\operatorname{corank}}
 \def\sign{\operatorname{sign}}
 \def\spn{\operatorname{span}}
 \def\R{{\mathbb R}}
 \def\Z{{\mathbb Z}}
 \def\0{{\mathbf 0}}
 \def\1{{\mathbf 1}}
 \def\a{{\mathbf a}}
 \def\b{{\mathbf b}}

\newcommand{\skcrossr}{\raisebox{-0.25\height}{\includegraphics[width=0.5cm]{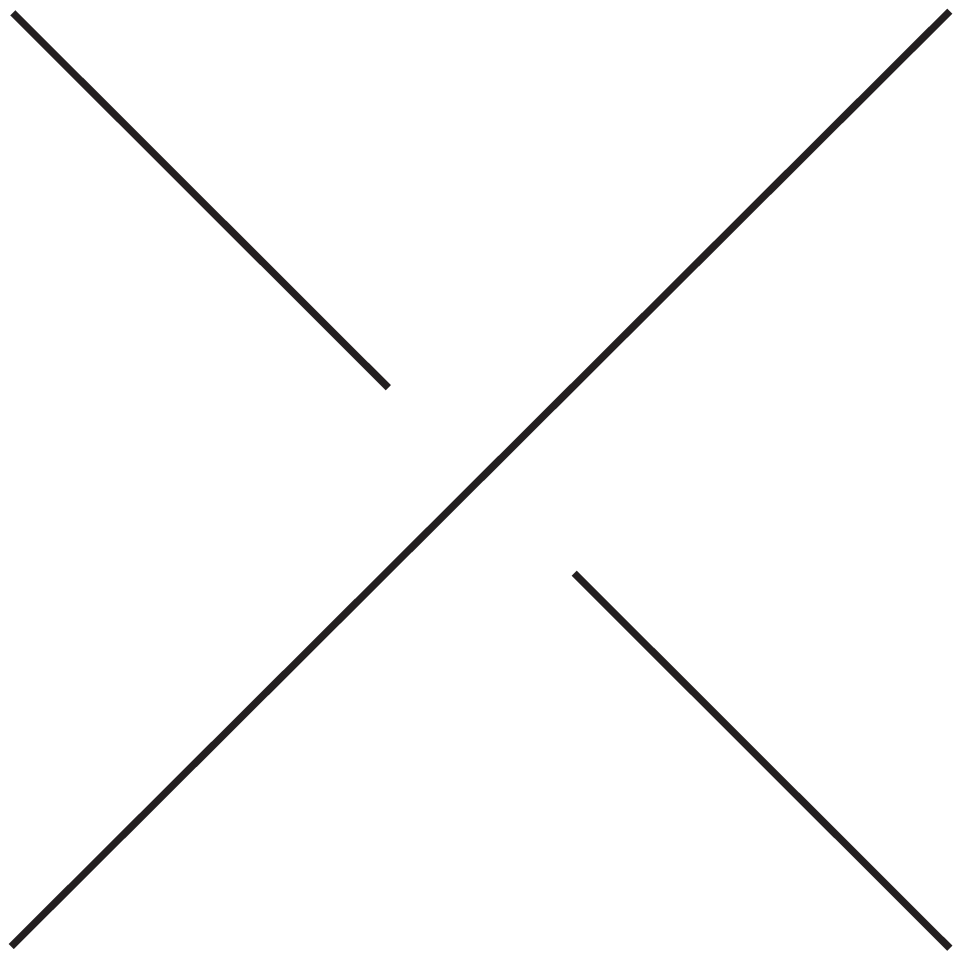}}}
\newcommand{\skcrv}{\raisebox{-0.25\height}{\includegraphics[width=0.5cm]{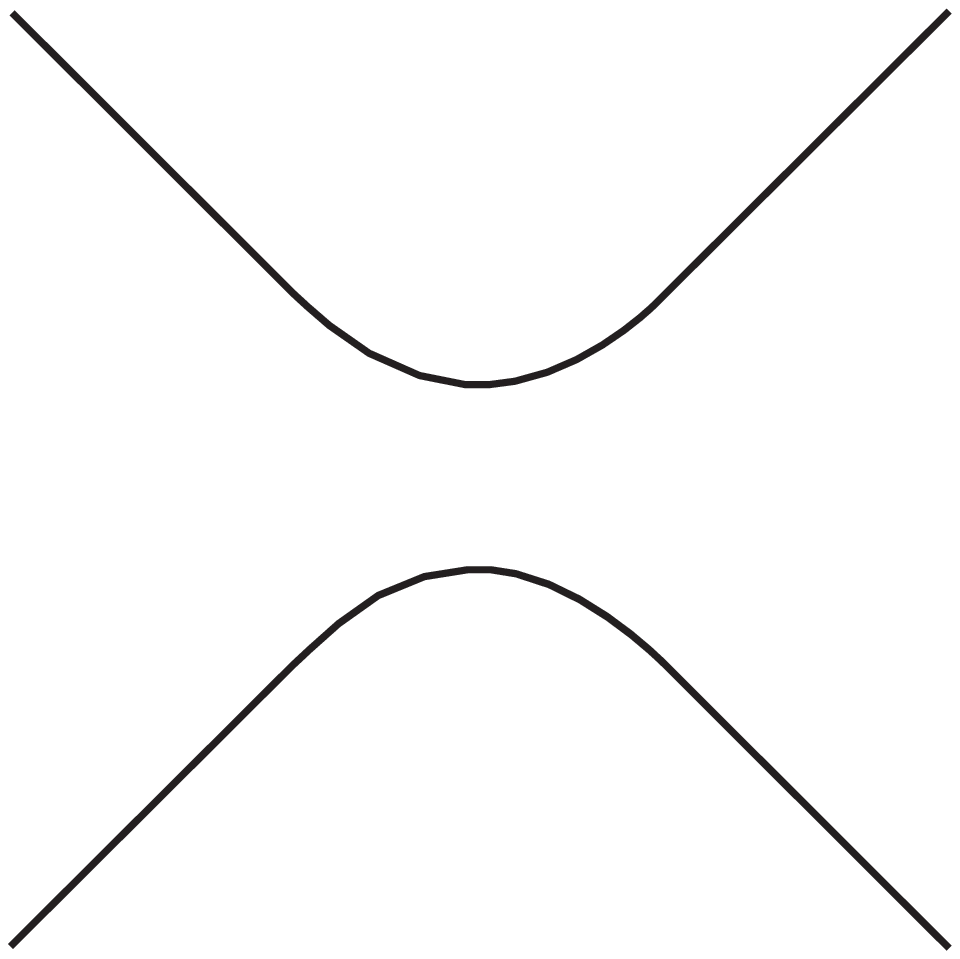}}}
\newcommand{\skcrh}{\raisebox{-0.25\height}{\includegraphics[width=0.5cm]{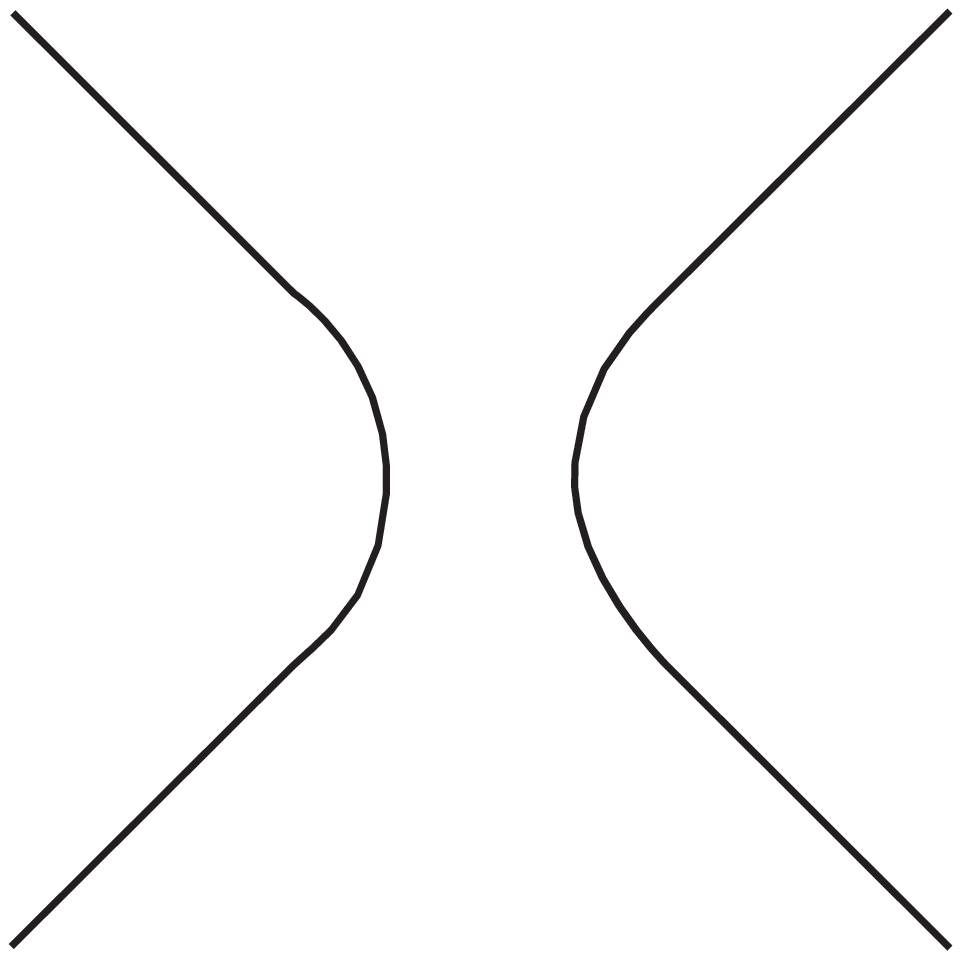}}}
\newcommand{\skcrro}{\raisebox{-0.25\height}{\includegraphics[width=0.5cm]{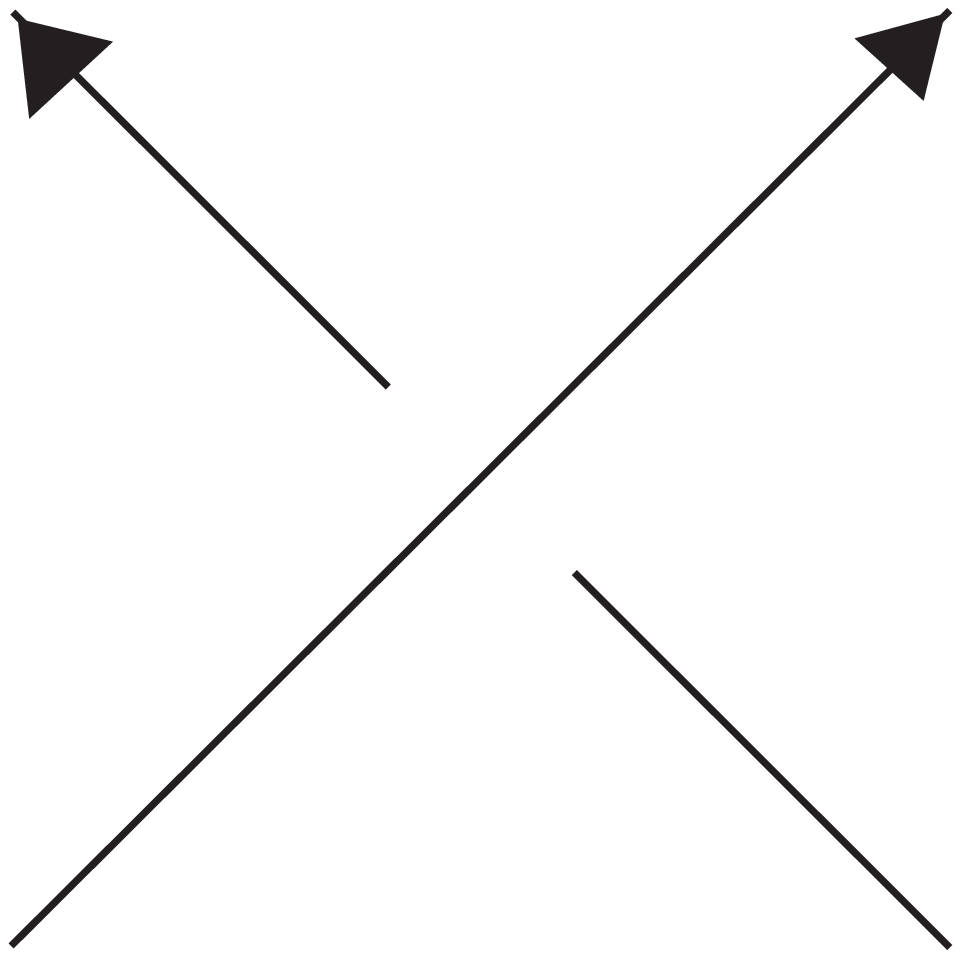}}}
\newcommand{\skcr}{\raisebox{-0.25\height}{\includegraphics[width=0.5cm]{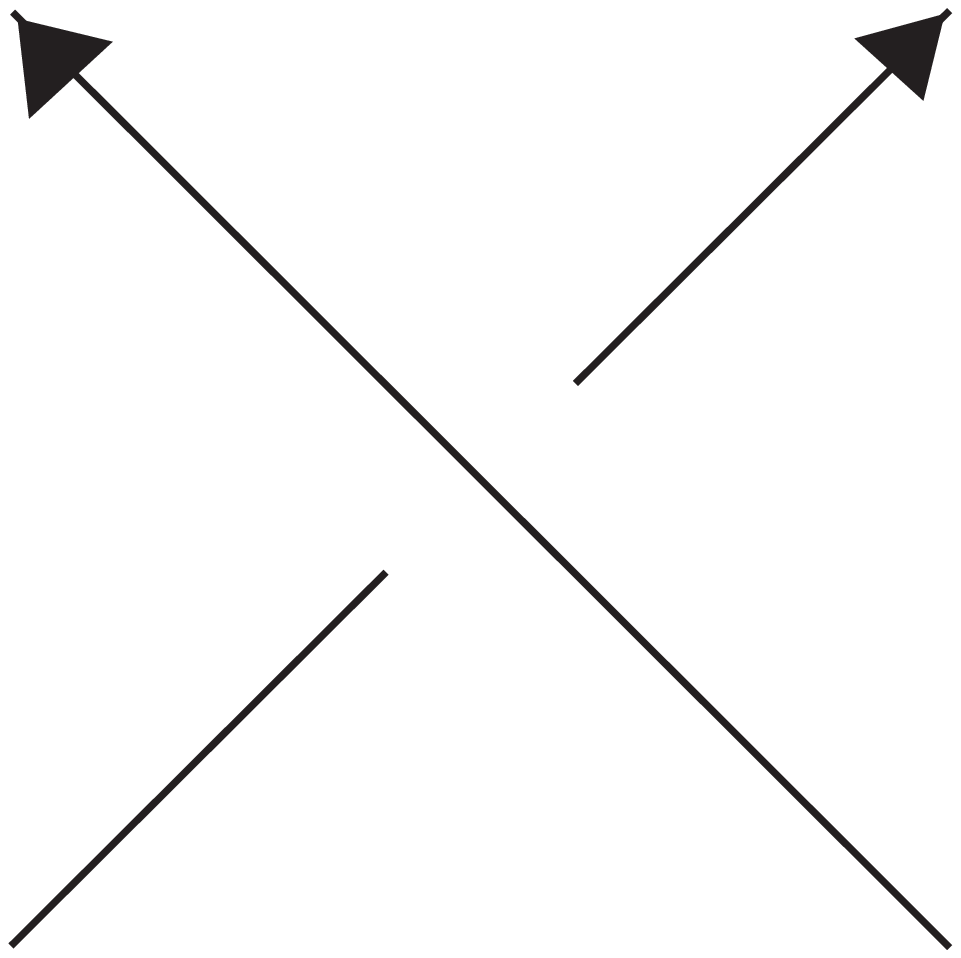}}}

\title {Introduction to Graph-Link Theory}

\author{D.\,P.~Ilyutko\footnote{Supported by grants of RF President
NSh -- 660.2008.1, RFFI 07--01--00648, RNP 2.1.1.7988.}\,,
V.\,O.~Manturov}
\date {}

\begin {document}

 \maketitle
 \abstract{The present paper is an introduction to a combinatorial theory
arising as a natural generalisation of classical and virtual knot
theory. There is a way to encode links by a class of `realisable'
graphs. When passing to generic graphs with the same equivalence
relations we get `graph-links'. On one hand graph-links generalise
the notion of virtual link, on the other hand they do not feel link
mutations. We define the Jones polynomial for graph-links and prove
its invariance. We also prove some a generalisation of the
Kauffman-Murasugi-Thistlethwaite theorem on `minmal diagrams' for
graph-links.}


 \section {Introduction}


The discovery of virtual knot theory by Kauffman~\cite{KaV} in
mid-1990-s was an important step in generalising combinatorial and
topological knot theoretical techniques into a larger domain (knots
in thickened surfaces), which is an important step towards
generalisation of these techniques for knots in arbitrary manifolds.

It turned out that some invariants (Kauffman bracket polynomial) can
be generalised for virtual knots immediately~\cite{KaV}, and some
other theories (Khovanov homology theory) need a complete revision
of the original construction for a generalisation for the case of
virtual knots~\cite{Izv}.

On the other hand, virtual knots sharpened several problems and
elicited some phenomena which do not appear in the classical knot
theory~\cite{FKM}, e.g., the existence of a virtual knot with
non-trivial Jones polynomial and trivial fundamental group
emphasises the difficulty of extracting the Jones polynomial
information out of the knot group.

In the present paper, we introduce a new class of objects closely
connected to both classical and virtual knots: graph-links. Likewise
virtual knots appear out of non-realisable Gauss code and thus
generalise classical knots (which have realisable Gauss codes),
graphs-links come out of {\em intersection graphs}: we may consider
graphs which realise chord diagrams, and, in turn, virtual links,
and pass to arbitrary simple graphs which correspond to some
mysterious objects generalising links and virtual links. Here we
refer to the paper by Traldi and Zulli~\cite{TZ} who made a similar
work but with another approach: there is a way of coding all knots
by {\em Gauss diagrams} and there is another way of coding all {\em
links} by {\em rotating circuits} (see, e.g.,~\cite{Arxiv}). Thus,
the class of objects described in the present paper is larger than
the one described by Traldi and Zulli. Rotating circuit approach to
both classical and virtual links has several advantages in
comparison to the Gauss diagram approach: for classical links, it
does not carry nugatory information, it is easier to recognise
planarity, and a single circle diagram can encode a link of
arbitrarily many components.

On the other hand, while passing from a chord diagram to its
intersection graph we forget a lot of information, say, mutant knots
are invisible at the level of intersection graphs. Thus, our graph
knots are a sort of `simplified generalised virtual knots'.

Though one can hardly imagine a picture of such a link for a
non-realisable graph, in the present paper we constructed an
invariant Kauffman bracket corresponding to them and Jones
polynomial. This bracket counts the ``number of non-existing circles
of a non-realisable chord diagram''. The invariance of the bracket
agrees with the fact that it does not feel knot mutation.

To construct the Kauffman bracket of a graph-link, one has to count
non-existing circles, but to manage with the Khovanov homology, one
should also know how these non-existing circles in different
Kauffman states of the diagram interfere. In general, the problem
seems to be unmanageable (e.g., since the Khovanov homology does
detect link mutation, Wehrli~\cite {Weh}), but the question is: to
which extent one can generalise (or categorify) the Kauffman bracket
of a knot.  We shall address this question in a sequel of the
present paper.

We also use various methods involving Kauffman bracket for detecting
the minimal crossing number (minimal vertex number) of a graph-knot,
and give many examples.


 \section {Preliminaries. Basic constructions}


  \subsection {Atoms and virtual links}

First of all, note that in this part of the paper we restrict
ourselves to a large class of graph-links corresponding to so-called
`orientable atoms'. Note that classical knots satisfy that
condition, i.e.\ they have an orientable atom. We shall deal with
the most general object in a sequel of the present paper.

 \begin {dfn}
An {\em atom} \cite {F} is a pair $(M,\Gamma)$ consisting of a
$2$-manifold $M$ and a graph $\Gamma$ embedded in $M$ together with
a colouring of $M\backslash \Gamma$ in a checkerboard manner. An
atom is called {\em orientable} if the surface $M$ is orientable.
Here $\Gamma$ is called the {\em frame} of the atom, whence by {\em
genus} (atoms and their genera were also studied by
Turaev~\cite{Turg}, and atom genus is also called the Turaev
genus~\cite {Turg}) ({\em Euler characteristic, orientation}) of the
atom we mean that of the surface $M$.
 \end {dfn}

 \begin {rk}
Throughout the paper, we shall deal with two types of graphs:
four-valent graphs (atom frames or shadows of links) and {\em
intersection-type graphs} of arbitrary valency. Four-valent graphs
are always assumed connected, and the intersection-type graphs may
not be connected, though they are assumed to have no loops and no
multiple edges.
\end {rk}

Having an atom, we may try to embed its frame in $\R^2$ in such a
way that the structure of opposite half-edges at vertices is
preserved. Then we can take the ``black angle'' structure of the
atom to restore the crossings on the plane (as ahead).

In~\cite{AtomsandKnots} it is proved that the link isotopy type does
not depend on the particular choice of embedding of the frame into
$\R^2$ with the structure of opposite edges preserved. The reason is
that such embeddings are quite rigid.

The atoms whose frame is embeddable in $\R^2$ with opposite
half-edge structure preserved are called {\em height} or {\em
vertical}.

However, not all atoms can be obtained from some classical knots.
Some abstract atoms may be quite complicated for its frame to be
embeddable into $\R^2$ with the opposite half-edges structure
preserved. However, if it is impossible to {\em embed} a graph in
$\R^2$, we may {\em immerse} it by marking artifacts of the
immersion (we assume the immersion to be generic) by small circles.

This leads to a connection between atoms and {\em virtual knots}
which perfectly agrees with {\em virtual knot theory} proposed by
Kauffman in~\cite{KaV}.

 \begin {dfn}
A {\em virtual diagram} is a $4$-valent diagram in $\R^2$ where each
crossing is either endowed with a classical crossing structure (with
a choice for underpass and overpass specified) or just said to be
virtual and marked by a circle.
 \end {dfn}

\begin {figure} \centering\includegraphics[width=250pt]{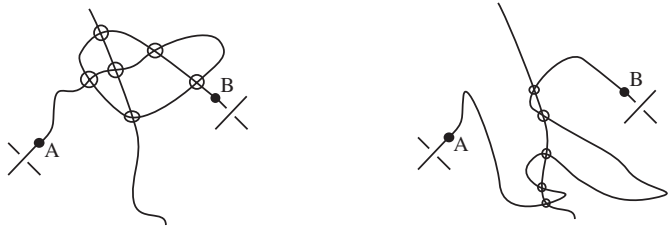}
  \caption{The detour move} \label{detour}
 \end {figure}

 \begin {dfn}
A {\em virtual link} is an equivalence class of virtual link
diagrams modulo generalised Reidemeister moves. The latter consist
of usual Reidemeister moves referring to classical crossings and the
{\em detour move} that replaces one arc containing only virtual
intersections and self-intersection by another arc of such sort in
any other place of the plane, see Fig.~\ref{detour}.
 \end {dfn}

 \begin {rk}
Throughout the paper, each virtual diagram has at least one
classical crossing.
 \end {rk}

 \begin {rk}
As the detour move does not affect the Kauffman bracket polynomial,
we make no difference between virtual diagrams obtained from each
other by detours. The coding by chord diagrams and graphs will not
see detours at all, so we will be able to check only classical
Reidemeister moves to establish any invariance.
 \end {rk}

\begin {figure} \centering\includegraphics[width=150pt]{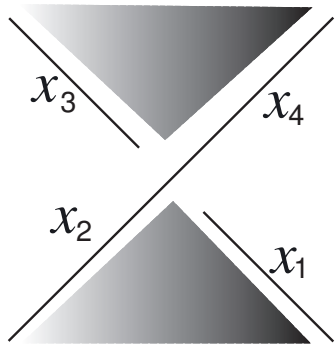}
  \caption{Gluing the cells} \label{black}
 \end {figure}

Given a virtual diagram $L$, let us construct the atom $A(L)$ out of
it. Vertices of $A(L)$ are in one-to-one correspondence with
classical crossings of the diagram $L$. Classical crossings of $L$
are connected to each other by branches of $L$ which may intersect
each other in virtual crossings. In each classical crossing we have
four emanating edges $x_1,x_2,x_3,x_4$ in clockwise-ordering such
that the pair $(x_1,x_3)$ forms an undercrossing and the pair
$(x_2,x_4)$ forms an overcrossing. These edges are in one-to-one
correspondence with edges of $A(L)$ connecting the corresponding
vertices. The $1$-cycles of the frame pasting black and white cells
are as follows. Each boundary of a $2$-cell is a rotating circuit on
a frame: a circuit which passes every edge at most once and switches
at each vertex from an edge to an adjacent (non-opposite) one. Black
cells are glued to the angles formed by $(x_1,x_2)$ and $(x_3,x_4)$,
and white cells are glued to the angles formed by $(x_2,x_3)$ and
$(x_1,x_4)$, see Fig.~\ref{black}. As a result we get an atom.

Let us consider the inverse operation. Having an atom, we may try to
construct a virtual knot diagram. To do that, we should take a
generic immersion of the atom's frame into $\R^2$, put virtual
crossings at the intersection points of images of different edges
and restore classical crossings at images of vertices `as above'.
Obviously, since we disregard virtual crossings, the most we can
expect is the well-definiteness up to detours. However, this allows
us to get different virtual link types from the same atom, since for
every vertex of the atom with four emanating half-edges $a,b,c,d$
(ordered cyclically on the atom) we may get two different
clockwise-orderings on the plane of embedding, $(a,b,c,d)$ and
$(a,d,c,b)$. This leads to a move called {\em virtualisation}.

\begin{figure} \centering\includegraphics[width=200pt]{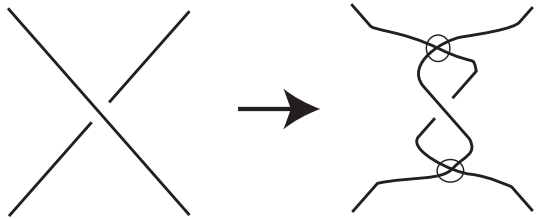}
  \caption{Virtualisation} \label{virtualisation}
 \end{figure}

 \begin {dfn}
By a {\em virtualisation} of a classical crossing of a virtual
diagram we mean a local transformation shown in
Fig.~\ref{virtualisation}.
 \end {dfn}

The above statements summarise as

 \begin{prop}\(see, e.g., {\em \cite{MBook}}\).
Let $L_1$ and $L_2$ be two virtual links obtained from the same atom
by using different immersions of its frame. Then $L_1$ differs from
$L_2$ by a sequence of \(detours and\/\) virtualisations.
 \end{prop}

Note that many famous invariants of classical and virtual knots
(Kauffman bracket, Khovanov homology~\cite{Kho,Izv},
Khovanov-Rozansky homology~\cite{KhR1,KhR2}) do not change under the
virtualisation, which supports the {\em virtualisation conjecture}:
if for two classical links $L$ and $L'$ there is a sequence
$L=L_{0}\to \cdots\to L_{n}=L'$ of virtualisations and generalised
Reidemeister moves then $L$ and $L'$ are classically equivalent
(isotopic).

Note that the usual virtual equivalence implies classical
equivalence for classical links, see~\cite{GPV}. This means that
classical links constitute a proper subset of virtual links.

\subsection{Thickened surface interpretation of virtual links}

Virtual links, being defined diagrammatically, have a topological
interpretation. They correspond to links in thickened surfaces
$S_{g}\times I$ with fixed $I$-bundle structure up to
stabilisations/destabilisations. Projecting  $S_{g}$ to $\R^{2}$
(with the condition, however, that all neighbourhoods of crossings
are projected with respect to the orientation, we get a diagram on
$\R^2$ from a generic diagram on $S_g$): besides the usual crossings
arising naturally as projections of classical crossings, we get
virtual crossings, which arise as artefacts of the projection: two
strands lie in different places on $S_g$ but they intersect on the
plane because they are forced to do so.

\begin{figure}
\centering\includegraphics[width=200pt]{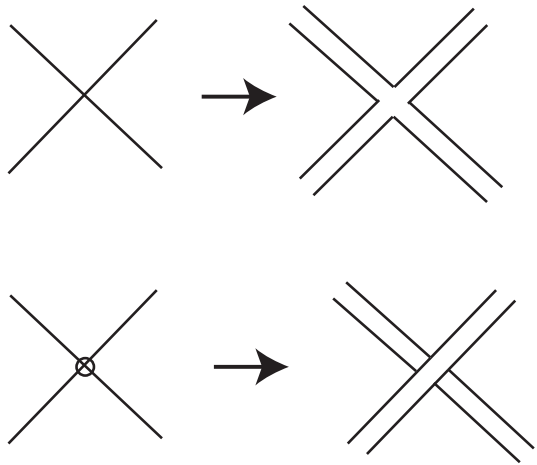} \caption{The
local structure}\label {crossskew}
 \end{figure}

We shall use the following construction. Having a virtual link
diagram $K$, we take all classical crossings of it and associate
with a neighbourhood of a crossing two crossing bands --- a `piece
of $2$-surface', as shown in Fig.~\ref{crossskew}, and with a
virtual crossing we associate a pair of skew bands, as shown in the
lower picture of Fig.~\ref{crossskew}.

If we connect these crossings and bands by (non-overtwisted) bands
going along edges, we get a $2$-surface with boundary. Gluing its
boundary components by discs, we get an orientable surface $M(K)$.
There is an obvious projection of (a part of) $M(K)$ to $\R^{2}$ in
such a way that classical crossings of the link diagram drawn in
$M(K)$ generate the diagram of $K$ (with virtual crossings at
intersections of skew bands). Then we get a link in $M(K)\times I$
representing the virtual link $K$

 \begin{rk}
Note that a thickened surface interpretation differs from the
presentation by atoms (see ahead): an atom need not be orientable,
while a checkerboard surface need not be $4$-colourable. On the
other hand, atom does not feel virtualisation, whence the thickened
surface does. Later, we will consider only connected thickened
surfaces and call the genus of it the {\em underlying genus} of a
representative for a virtual knot $K$.
 \end{rk}

For more details, see, e.g.,~\cite {MBook}.

 \subsection {Chord diagrams and intersection graphs}

  \begin{dfn}
A {\em chord diagram} is a graph consisting of a selected oriented
cycle (the {\em circle}) and several non-oriented edges ({\em
chords}) connecting points on the circle in such a way that every
point is incident to at most one chord. A chord diagram is {\em
labeled} if every chord is endowed with a sign `$+$' or `$-$'. If no
signs are indicated, we assume the chord diagram has all chords with
sign `$+$'. A labeled chord diagram is a {\em $d$-diagram} if the
corresponding intersection graph is bipartite. Two chords in a chord
diagram are called {\em linked} if the ends of one chord lie in
different connected components of the circle without the end-points
of the second chord.
  \end{dfn}

Some remarks are in order.

 \begin{rk}
Note that $d$-diagrams are precisely those encoding classical link
diagrams \cite{AtomsandKnots}: a chord diagram is embeddable in
$\R^2$ iff it is a $d$-diagram, and embeddability of a chord diagram
yields the planarity of the corresponding four-valent graph (a
shadow of the link, which thus has no virtual crossings).
 \end{rk}

  \begin{rk}
To avoid confusion between framing and labeling throughout the
paper, we consider only labeled chord diagram with positive framing
as each virtual link corresponds to an orientable atom and we forget
about framing.
  \end{rk}

 \begin{rk}
In fact, chord diagrams with all positive chords encode all
orientable atoms with one white cell: this white cell corresponds to
the only $A$-state of the virtual diagram, and chords show how this
cell approaches itself in neighbourhoods of crossings (atom
vertices). If we want to deal with all atoms and restrict ourselves
for the case of one circle, we should take this circle to be
corresponding to some other state of the atom, which is encoded by
labelings of the chords.
 \end{rk}

Having a labeled chord diagram $D$, one can construct a virtual link
diagram $K(D)$ (up to virtualisation) as follows. Let us immerse
this diagram in $\R^{2}$ taking an embedding of the circle and
placing some chords inside the circle and the other ones outside the
circle. After that we remove a neighbourhoods of each of the chord
ends and replace it by a couple of lines with a classical crossing
in the following way. A crossing can be smoothed in two ways: $A$
and $B$ as in the Kauffman bracket polynomial. We require that the
initial piece of the circle corresponds to the $A$-smoothing if the
chord is positive and to the $B$-smoothing if it is negative:
$A:\skcrossr\to\skcrh$, $B\:\skcrossr\to\skcrv$. Note that from
$d$-diagrams we get classical links.

Conversely, having a connected virtual diagram $K$ (having at least
one classical crossing) with oriented atom, one gets a labeled chord
diagram. Indeed, one takes a circuit of $K$ which is a map from
$S^{1}$ to the shadow of $K$ which is bijective outside classical
and virtual crossings having exactly two inverse images at each
classical crossing and two virtual crossings, going transversally at
each virtual crossings, and turning from an edge to an adjacent
(non-opposite) edge at each classical crossing. This defines a chord
diagram (orientation of the circle is chosen arbitrarily), where the
sign of the chord is `$+$' if the circuit locally agrees with the
$A$-smoothing, and `$-$' if it agrees with the $B$-smoothing.

It can be easily checked that this operation is indeed inverse to
the operation of constructing a virtual link out of a chord diagram:
if we take a chord diagram $D$, and construct a virtual diagram
$K(D)$ out of it, then for some circuit the chord diagram
corresponding to $K(D)$ will coincide with $D$. The rule for setting
classical crossings here agrees with the rule described above
because of the orientability of the atom.

This proves the following

  \begin {theorem}[\cite {MBook}]
Any connected virtual diagram with an orientable atom is obtained
from a certain labeled chord diagram.
  \end {theorem}

We restrict our attention to connected virtual diagrams with an
orientable atoms. The following two theorems show that the set of
all such diagrams gives us another class which is wider than the
class of classical links but narrower than the one of all virtual
links.

As the goal of our paper is to construct the Kauffman bracket
polynomial and the Jones polynomial, we disregard virtualisation.

 \begin{theorem}
Each two equivalent \(in the class of all virtual diagrams\/\)
connected virtual diagrams are equivalent in the class of connected
virtual diagrams.
 \end{theorem}

 \begin{sproof}
Having a disconnected diagram, we can make it connected by adding
crossings by second Reidemeister moves. After performing a necessary
move, we may perform the inverse Reidemeister move (and, possibly,
perform another move in another place). It is easy to see that the
process can be organised in such a way that the ``auxiliary''
crossings do not affect the original moves we have to perform, and
in the final diagram all auxiliary crossings are removed. This
completes the proof.~$\square$
 \end{sproof}

 \begin{theorem}
Each two equivalent \(in the class of all virtual diagrams\/\)
virtual diagrams $K_1$ and $K_2$ with orientable atoms are
equivalent in the class of virtual diagrams with orientable atoms.
\end{theorem}

This theorem was independently proved by O.\,Ya.~Viro and the second
named author of the present paper, but the proof was not published.
This theorem signifies the consistency of the class of objects we
are dealing with: virtual links with orientable atoms constitute a
proper part of all virtual links.

 \begin{proof}
Consider a thickened surface presentation of virtual links,
\cite{KaV,Kup}. It can be easily checked that the orientability of
the atom corresponding to a virtual link diagram depends only on the
$\Z_{2}$-homology class of the corresponding link surface
presentation. Consequently, if we perform a move which does not
change the genus of the corresponding surface, we do not change
orientability of the atom. Obviously, the detours do not affect the
surface presentation and hence orientability of the atom. A
straightforward check shows neither the first nor the third
Reidemeister move changes the genus of the underlying surface, thus,
does affect the orientability of the atom. The second Reidemeister
move may either increase the genus or decrease the genus or leave it
as it is. For the case of the second Reidemeister move not changing
the genus, the orientability is not changed. Moreover, we may pass
from an orientable atom to an orientable atom only in the case when
we increase the genus of the thickened surface representation by a
second Reidemeister move. Now, assume $K_{1}$ and $K_{2}$ are
virtual diagrams with orientable atoms. We use the same notation for
the corresponding knots in thickened surfaces: $K_{1}$ and $K_{2}$.
It follows from Kuperberg's theorem that there is the following
sequence of transformations: we first transform $K_{1}$ to its
minimal representative $K'_{1}$ having minimal genus. Then we
transform $K_{2}$ to its minimal representative $K'_{2}$ which has
the same genus as $K'_1$ and is isotopic inside the surface of that
genus to $K'_{1}$. The latter means that $K'_{1}$ and $K'_{2}$ are
connected by a sequence of genus-preserving Reidemeister moves. So,
since both $K_{1}$ and $K_{2}$ have orientable atoms, the whole
chain of diagrams $K_{1}\to\dots\to K'_{1}\to\dots\to
K'_{2}\to\dots\to K_{2}$ corresponds to a sequence of orientable
atoms, which completes the proof.
 \end{proof}

The Reidemeister moves give a combinatorial description of the
relationship between the different diagrams of a given virtual link.
The following definition describes the set of moves corresponding to
usual set of Reidemeister moves for planar diagrams. After that we
next pass on to moves on graphs.

\begin{figure} \centering\includegraphics[width=250pt]{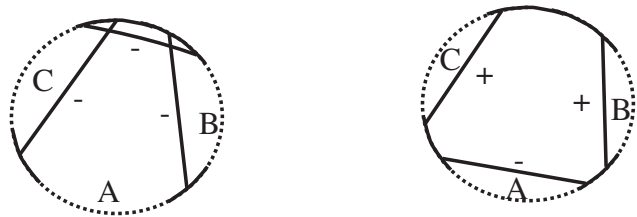}
  \caption{The third Reidemeister chord-move} \label{third move}
 \end{figure}

 \begin{dfn}
$\Omega 1$. The first Reidemeister chord-move is an addition/removal
of an isolated chord labeled `$+$' or `$-$' (an isolated chord, not
linked with any others).

$\Omega 2$. The second Reidemeister chord-move is an
addition/removal of two parallel chords labeled `$+$' and `$-$', so
that these chords have the same linked with others.

$\Omega 3$. The third Reidemeister chord-move is shown in Fig.~\ref
{third move}. All the chords except three chords shown in  Fig.~\ref
{third move} are fixed and their ends lie on the `dotted' parts of
the circle.

$\Omega 4$. The fourth chord-move is shown in Fig.~\ref {fourth
move}. The move takes four segments of the chord diagram denoted by
$A,B,C,D$ and  transforms the diagram as follows. We perform the
surgery along the chords with signs $a$ and $b$, $a,b\in\{\pm1\}$.
This surgery cuts the diagram into four pieces $A,B,C,D$ each
containing some chord ends, and reconnects them by a pair of
vertical lines and a pair of horizontal lines (middle picture). The
two chords we are operating with change their labels. After that we
redraw the figure to get a `round circle' of the chord diagram
(rightmost picture) and get a shuffle of the segments $A,B,C,D$ with
all chord ends lying on them.
 \end{dfn}

\begin{figure} \centering\includegraphics[width=250pt]{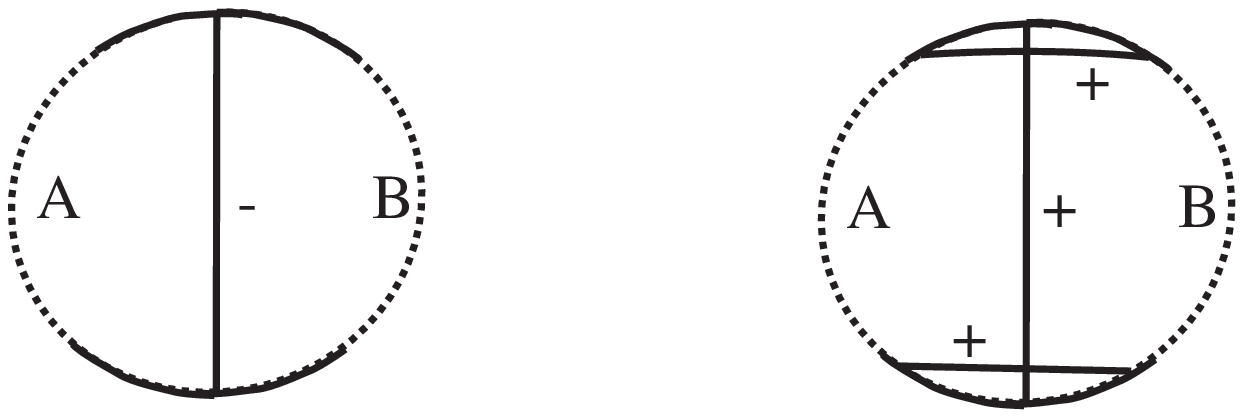}
  \caption{The fifth chord-move} \label{fifth move}
 \end{figure}

 \begin{rk}
Using the second Reidemeister and fourth chord-moves, we can replace
each chord labeled `$-$' by three chords labeled `$+$' as shown in
Fig.~\ref {fifth move}.  Let us call this transformation the {\em
fifth move} and denote it by $\Omega 5$. The chords having their
ends on the `dotted' parts of the circle are fixed.
 \end{rk}

 \begin{figure} \centering\includegraphics[width=320pt]{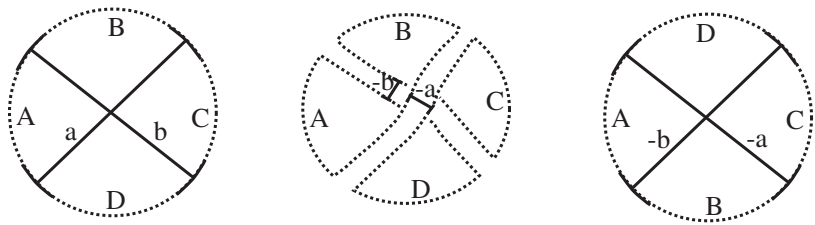}
  \caption{The fourth chord-move} \label{fourth move}
 \end{figure}

 \begin {theorem}\label {th:link_chord}
Let $K_1$ and $K_2$ be two connected virtual diagrams with
orientable atoms, and $D_1$ and $D_2$ be two labeled chord diagrams
obtained from $K_1$ and $K_2$, respectively. If $K_1$ and $K_2$ are
equivalent in the class of connected diagrams with orientable atoms
then $D_1$ is obtained from $D_2$ by $\Omega 1 - \Omega 4$
chord-moves.
 \end {theorem}

\begin {figure} \centering\includegraphics[width=250pt]{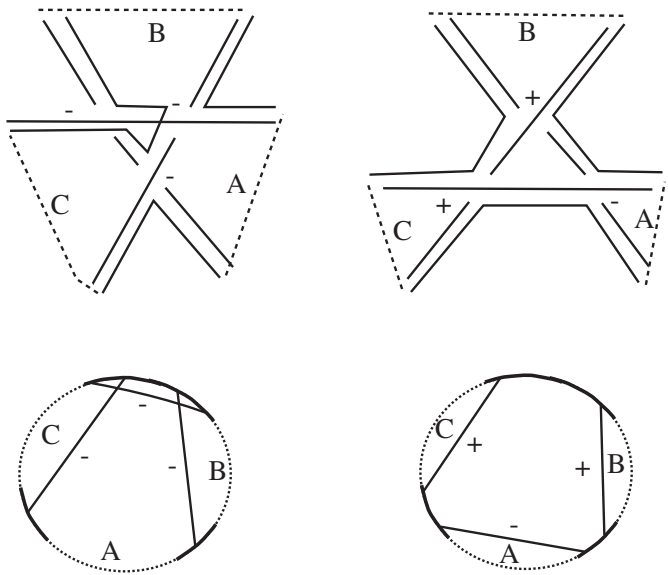}
  \caption{The third Reidemeister move and the third chord-move}\label{3rm move}
 \end {figure}

 \begin {proof}
First of all, note that the chord diagrams of two virtual links
$K_1$ and $K_2$ distinguished from each other by the detour moves
are the same.

Secondly, when we construct a chord diagram from a virtual diagram
$K$ we fix a circuit of $K$ which is a map from $S^{1}$ to the
shadow of $K$. The independence from the choice of a circuit is
guaranteed by the fourth chord-move $\Omega 4$. Indeed, assume two
labeled chord diagrams $D_{1}$ and $D_{2}$ represent the same
virtual diagram (up to detours). Then $D_{1}$ and $D_{2}$ can be
obtained from each other by a change of rotating circuit. We have to
rotate the circuit for $D_{1}$ at some collection of $m$ vertices
(chords of $D_{1}$) to get the circuit for $D_{2}$. It is easy to
see that there is a couple of linked chords among them (otherwise
there is no circuit, but a collection of circuits after a surgery).
Take a pair of intersecting chords and change the circuit
respectively. This will be precisely one fourth move. We get a
diagram $D'_{1}$ which differs from $D_{2}$ at some collection of
$m-2$ vertices. Reiterating the process above, we get to $D_{2}$
from $D_{1}$ in a finite number of steps.

Without loss of generality, we can assume that $K_1$ and $K_2$ are
distinguished from each other by classical Reidemeister moves, for
each $K_1$ and $K_2$ we have rotating circuits of $K_1$ and $K_2$,
and all the moves are performed in the class of connected virtual
diagrams with orientable atoms.

Under such assumptions the classical first and second Reidemeister
moves on virtual link are equivalent to the first and second
Reidemeister chord-moves.

As for the third Reidemeister move on virtual diagrams, it is
sufficient to consider only one variant provided that we have the
whole collection of the first and second Reidemeister moves,
see~\cite{Oht}. At the level of chord diagrams, it is again
sufficient to consider only one way of representing this third
Reidemeister move (shown in Fig.~\ref{3rm move}) provided that we
have all second Reidemeister chord-moves and the fourth chord-move.
Under such assumptions the one classical third Reidemeister move on
virtual link is equivalent to the third Reidemeister chord-move,
Fig.~\ref {3rm move}. This completes the proof.
 \end {proof}

Assume we are given a labeled chord diagram $D$. Let us construct
the {\em labeled intersection graph}, see~\cite{CDL}, $G(D)$ as
follows. This is the labeled graph whose vertices are in one-to-one
correspondence with chords of the diagram, the label of each vertex
corresponding to a chord coincides with that of the chord, and two
vertices are connected by an edge if and only if the corresponding
chords are linked.

 \begin {dfn}
Let $G$ be a labeled graph on $n$ vertices.  Fix an enumeration of
vertices for $G$. We define the {\em adjacency matrix} $A(G)$ as
follows. Set $A(G)=(a_{ij})_{i,j=1,\dots,n}$ be the {\em adjacency
matrix} of $G$, i.e.\ $a_{ij}=1$ if and only if the vertices $i$ and
$j$ are incident, and $a_{ij}=0$ otherwise. Besides, we
set\footnote{The case $a_{ii}=1$ corresponds to framed chords with
framing $1$, which, in turn, correspond to non-orientable atoms,
that we shall consider in another paper.} $a_{ii}=0$.
 \end {dfn}

\begin {figure} \centering\includegraphics[width=200pt]{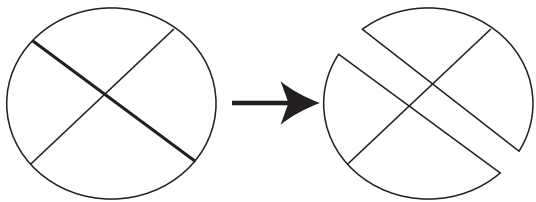}
  \caption{A surgery of the circuit along a chord} \label{surger}
 \end {figure}

\begin {figure} \centering\includegraphics[width=230pt]{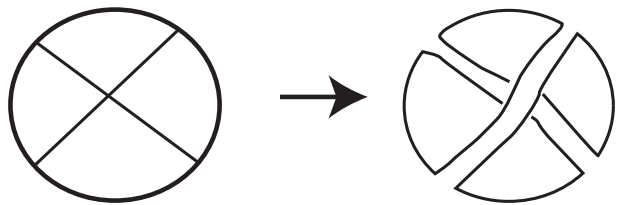}
  \caption{The manifold $m(D)$} \label{result}
 \end {figure}

Assume we are given a labeled chord diagram. Define the {\em surgery
over the set of chords as follows}. For every chord, we draw a
parallel chord near it and remove the arc of the circle between
adjacent ends of the chords as in Fig.~\ref {surger}. By a small
perturbation, the picture in $\R^{2}$ is transformed into
one-manifold in $\R^{3}$. This manifold $m(D)$ is the result of
surgery, see Fig.~\ref {result}.

Surprisingly, the number of the connected components of $m(D)$ can
be counted out of the intersection graph.

  \begin {theorem}[\cite {Soboleva},~\cite {BNG}]\label {th:sob}
Let $D$ be a \(labeled\/\) chord diagram, and let $G$ be its
\(labeled\/\) intersection graph. Then the number of connected
components of $m(D)$ equals $\corank A(G)+1$, where $A(G)$ is the
adjacency matrix of $G$.
  \end {theorem}

 \begin {rk}
Theorems~\ref {th:sob} allows us to define `the number of circles'
for a graph even when the given graph is not a intersection graph.
 \end {rk}

\subsection{Acknowledgments}

The authors are grateful to L.\,H.~Kauffman, V.\,A.~Vassiliev, and
A.\,T.~Fomenko for their interest to this work.


 \section {The Kauffman bracket polynomial}


The Kauffman bracket polynomial~\cite{KaV,KauffmanBracket} is a very
useful model for understanding the Jones polynomial~\cite{Jones}.
The Kauffman bracket polynomial associates with a virtual link
diagram a Laurent polynomial in one variable $a$. After a small
normalisation (multiplication by a power of $(-a)$) it gives an
invariant for virtual links.

\begin {figure}
  \centering\includegraphics[width=200pt]{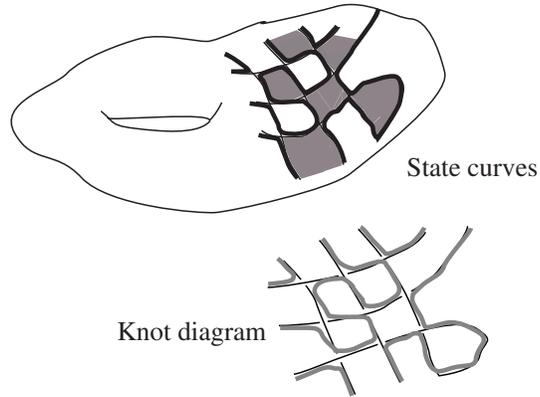} \caption{State
  curves drawn on an atom} \label{astateonanatom}
 \end {figure}

This invariant can be read out of the atom corresponding to a knot
diagram. Namely, take an atom $V$ with $n$ vertices corresponding to
a virtual diagram $L$ with $n$ classical crossings. By a {\em state}
we mean a choice of a pair of black or white angles at every vertex
of $V$. Every such choice gives rise to a collection of closed
curves on $V$ whose boundaries contain all  edges of $V,$ see
Fig.~\ref{astateonanatom}, and at each crossing the curves turn
locally from one edge to an adjacent edge sharing the same angle of
the prefixed colour.

Thus, having $2^{n}$ states of the atom, we define the Kauffman
bracket polynomial of it as
 \begin {equation}
\langle V\rangle
=\sum_{s}a^{\alpha(s)-\beta(s)}(-a^{2}-a^{-2})^{\gamma(s)-1},\label{Kauffmanbracket}
 \end {equation}
where the sum is taken over all states $s$ of the diagram,
$\alpha(s)$ and $\beta(s)$ denote the number of white and black
angles in the state (thus, $\alpha(s)+\beta(s)=n$ and $\gamma(s)$
denotes the number of curves in the state).

As mentioned above, the Kauffman bracket polynomial is invariant
under the virtualisation. Thus, it is not surprising that it can be
read out of the corresponding atom.

If the atom $A$ is obtained from a (framed) chord diagram $C$, then
one can construct the Kauffman bracket polynomial $\langle
C(A)\rangle$.

Thus, one obtains a function $f$ on (framed) chord diagrams valued
in Laurent polynomials in $a$. We shall return to that function
because it is connected to the Vassiliev invariants of knots and
$J$-invariants of closed curves (Lando,~\cite{Lando}).

Throughout the paper, we consider `intersection-type' graphs without
loops and multiple edges ({\em strict} graphs). Assume now we have a
labeled graph (a graph with each vertex labeled either positively or
negatively).

\begin {figure}
  \centering\includegraphics[width=200pt]{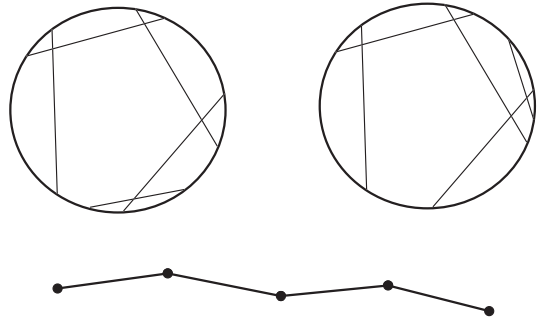} \caption{Two
  chord diagrams with the same intersection graph}
  \label{nonuniqueness}
 \end {figure}

Then it may or may not be represented as an intersection graph of a
chord diagram (see~\cite{Bouchet} for the details) for which it is
an intersection graph. Moreover, if such a chord diagram exists, it
should not be unique, see, e.g., Fig.~\ref {nonuniqueness}.

This non-uniqueness usually corresponds to so called {\em mutations}
of virtual knots.

\begin {figure} \centering\includegraphics[width=200pt]{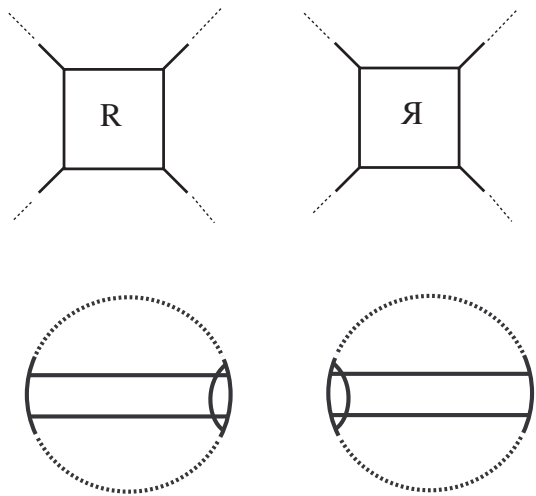}
  \caption{Mutation and its chord diagram presentation}\label{muta}
 \end {figure}

The mutation operation (shown in the top of Fig.~\ref{muta}) cuts a
piece of a knot diagram inside a box, turns it by a half-twist and
returns to the initial position. Note that the virtualisation is a
special case of the mutation.

It turns out that the mutation operation is expressed in terms of
chord diagrams in almost the same way: one cuts a piece of diagram
with $4$ ends and exchanges the top and the bottom parts of it (see
bottom picture of Fig.~\ref{muta}). Exactly this operation
corresponds to the mutation from both Gauss diagram and rotating
circuit points of view.

In the bottom part of Fig.~\ref{muta} chords whose end points belong
to the `dotted' area remain the same; the other chords (having their
end-points in the `dotted' segments) are reflected as a whole.

It is well-known that the Kauffman bracket polynomial does not
detect mutations, see, e.g.,~\cite {MBook}. Thus, one might guess
that the corresponding Kauffman bracket polynomial can be read out
of the intersection graph.

Surprisingly, the Kauffman bracket polynomial can be defined (the
initial definition was given in~\cite {Arxiv}, but graph-links and
Reidemeister moves for graphs were not defined and there was no
approach to prove any invariance) in a meaningful way even for those
labeled graphs which can not be represented as intersection graphs
of chord diagrams.

As we are able to calculate how many circles we have in each state,
we can apply  (\ref{Kauffmanbracket}) to calculate the Kauffman
bracket polynomial.

 \begin {dfn}
Let $G$ be a (strict) labeled finite graph with the set of vertices
$V(G)$ and the set of edges $E(G)$. Suppose $s\subset V(G)$. Set
$G(s)$ to be the subgraph of the graph $G$ with the set of vertices
$V(G(s))=s$ and the set of edges $E(G(s))$ such that $\{u,v\}\in
E(G(s))$, where $u,v\in s$, iff $\{u,v\}\in E(G)$.
 \end {dfn}

 \begin {dfn}
We call a subset of $V(G)$ a {\em state} of the graph $G$. A state
is called the $A$-{\em state} if it consists of all the vertices of
$G$ labeled `$-$' and no vertex labeled `$+$'. Analogously, a state
is called the $B$-{\em state} if it consists of all vertices of $G$
labeled `$+$' and no vertex labeled `$-$'. The {\em opposite} state
for a state $s$ is the set of vertices complement to $s$ (the
opposite state to the $A$-state is the $B$-state). Two states are
called {\em neighbouring} if they differ only in one vertex, which
belongs to one state and no to other. The {\em distance} between two
states is equal to the number of the vertices in which two states
differ. A state $s$ is called {\em single-circle} state if $\corank
A(G(s))=0$. We define the {\em number of circles} in a state $s$ as
$\corank A(G(s))+1$.
 \end {dfn}

 \begin {rk}
Note that $k+l\leqslant n+2$, here $k$ and $l$ are the numbers of
circles in $A$-state and $B$-state, respectively.
 \end {rk}

 \begin {dfn}
Let $G$ be a strict finite labeled graph with the set of vertices
$V(G)$, $|V(G)|=n$. Define the Kauffman bracket polynomial of $G$ as
  $$
\langle
G\rangle=\sum\limits_{s}a^{\alpha(s)-\beta(s)}(-a^2-a^{-2})^{\corank
A(G(s))},
  $$
where the sum is taken over all states $s$ of the graph $G$,
$\alpha(s)$ is equal to the sum of the vertices labeled `$-$' from
$s$ and of the vertices labeled `$+$' from $V(G)\setminus s$,
$\beta(s)=n-\alpha(s)$. Recall that all matrices are over $\Z_2$ and
$\corank A_{n\times n}=n-\operatorname{rank}A$.
 \end {dfn}


 \section {Graph-links}


One can construct the intersection graph from a chord diagram.
Therefore we can define graph-moves corresponding the chord-moves.
As a result we have a new object --- an equivalence class of labeled
graphs over formal moves.

 \begin {dfn}
Let $G$ be a graph and let $v$ be its vertex. The set of all
vertices adjacent to $v$ is called {\em environment of a vertex} and
denoted by $N(v)$ or $N_G(v)$.
 \end {dfn}

 \begin {dfn}
$\Omega_g 1$. The first Reidemeister graph-move is an
addition/removal of an isolated vertex labeled `$+$' or `$-$'.

$\Omega_g 2$. The second Reidemeister graph-move is an
addition/removal of two vertices not connected by an edge and having
the different signs and the same adjacency with others.

$\Omega_g 3$. The third Reidemeister graph-move is defined as
follows. Let $u,\,v,\,w$ be three vertices of $G$ all having label
`$-$' and $u$ is adjacent in $G$ only to $v$ and $w$. Then we change
only the adjacency of $u$ with the vertices $t\in N(v)\setminus
N(w)\bigcup N(w)\setminus N(v)$ (for other pairs of vertices we do
not change their adjacency). In addition, we change the labels of
$v$ and $w$ to `$+$'. The inverse operation is also called the third
Reidemeister graph-move.

$\Omega_g 4$. The fourth graph-move is defined as follows. We take
two adjacent vertices $u$ labeled $a$ and $v$ labeled $b$,
$a,\,b\in\{\pm1\}$. Then we change labels of $u$ and $v$ so that the
label of $u$ becomes $-b$ and the label of $v$ becomes $-a$. In
addition, we change the adjacency for each pair $(t,w)$ of vertices,
where $t\in N(u),\,w\in N(v)\setminus N(u)$ and $t\in N(v),\,w\in
N(u)\setminus N(v)$.
 \end {dfn}

 \begin {rk}
The fifth graph-move analogous to the fifth chord-move is as
follows. We take a vertex $u$ labeled `$-$', change its label and
add two vertices labeled `$+$' and joined with $u$. This operation
is denoted by $\Omega_g 5$.
 \end {rk}

The definition above together with the arguments of Theorem~\ref
{th:link_chord} yields

 \begin {theorem}\label {th:link_graph}
Let $K_1$ and $K_2$ be two connected virtual diagrams with
orientable atoms, and $G_1$ and $G_2$ be two labeled intersection
graphs obtained from $K_1$ and $K_2$, respectively. If $K_1$ and
$K_2$ are equivalent in the class of connected diagrams with
orientable atoms then $G_1$ is obtained from $G_2$ by $\Omega_g 1 -
\Omega_g 4$ graph-moves.
 \end {theorem}

We defined the Kauffman bracket polynomial. Our goal is to show that
it is invariant under some graph-moves, and then to normalise it in
order to get the Jones polynomial. As a result, we get an invariant
of graph-links, see ahead. Now we are going to define this object
and to check the invariance of the Kauffman bracket polynomial up to
some factor.

 \begin {dfn}
A {\em graph-link} is an equivalence of strict labeled graphs modulo
$\Omega_g 1 - \Omega_g 4$ graph-moves.
 \end {dfn}

 \begin {theorem}\label {th:mov_reid}
The Kauffman bracket polynomial of a labeled graph is invariant
under $\Omega_g 2 - \Omega_g 4$ graph-moves.
 \end {theorem}

 \begin {proof}
Let $G$ be a labeled graph with the set of vertices $V(G)$ and let
$\widetilde{G}$ be a graph obtained from $G$ by some graph-move of
the graph-moves $\Omega_g 2-\Omega_g 4$. $V(\widetilde{G})$ is the
set of all vertices of $\widetilde{G}$.

  \begin {enumerate}
   \item
Let us consider the second graph-move $\Omega_g 2$. Suppose that
vertices $u\in V(\widetilde{G})$ labeled `$+$' and $v\in
V(\widetilde{G})$ labeled `$-$' are added to $G$ to get
$\widetilde{G}$. Enumerate all vertices of $G$ by numbers from $1$
to $n$, and enumerate all vertices of $\widetilde G$ by numbers from
$1$ to $n+2$ so that the number of $u$ is $1$, the number of $v$ is
$2$, and $w\in V(\widetilde{G})\setminus\{u,v\}$ has number
$i\geqslant 3$ if and only if $w$ in $G$ has number $(i-2)$. The
adjacency matrix $A(\widetilde{G})$ is
  $$
A(\widetilde{G})=\left(\begin {array} {ccc}
 0& 0 & \a^\top\\
 0& 0 & \a^\top\\
 \a& \a& A(G)
 \end {array}\right),
  $$
where bold characters indicate row and column vectors.

Let $s\subseteq V(G)$. Then
 \begin {gather*}
\corank A(\widetilde G(s))=\corank A(G(s)),\\
\corank A(\widetilde G(s\cup\{u\}))=\corank A(\widetilde
G(s\cup\{v\}))=\corank A(\widetilde G(s\cup\{u,v\}))-1
 \end {gather*}
and
   \begin {gather*}
\langle\widetilde{G}\rangle=\sum\limits_{s\subseteq
V(G)}\Bigl(a^{2(\alpha(s)+1)-(n+2)}(-a^2-a^{-2})^{\corank A(
\widetilde G(s))}+\\
+a^{2\alpha(s)-(n+2)}(-a^2-a^{-2})^{\corank A(\widetilde
G(s\cup\{u\}))}+\\
+a^{2(\alpha(s)+2)-(n+2)}(-a^2-a^{-2})^{\corank A(\widetilde
G(s\cup\{v\}))}+\\
+a^{2(\alpha(s)+1)-(n+2)}(-a^2-a^{-2})^{\corank A(\widetilde
G(s\cup\{u,v\}))}\Bigr)=\\
=\sum\limits_{s\subseteq
V(G)}\Bigl(a^{2\alpha(s)-n}(-a^2-a^{-2})^{\corank
A(G(s))}+\\
+a^{2\alpha(s)-n}(-a^2-a^{-2})^{\corank A(\widetilde
G(s\cup\{u\}))}(a^{-2}+a^2-a^2-a^{-2})\Bigr)=\langle G\rangle,
   \end {gather*}
here $\alpha(s)$ is equal to the sum of the vertices labeled `$-$'
from $s$ and of the vertices labeled `$+$' from $V(G)\setminus s$.

  \item
Let us consider the third graph-move $\Omega_g 3$. There are two
versions of the third graph-move we consider only one. Enumerate all
vertices of $G$ by numbers from $1$ to $n$ so that the label of $u$
is $1$, the label of $v$ is $2$, and the label of $w$ is $3$. The
vertices of $\widetilde{G}$ have the same numbers as the vertices of
$G$ do. The adjacency matrices $A(G)$ and $A(\widetilde{G})$ are
  $$
A(G)=\left(\begin {array} {cccc}
 0 & 1 & 1 & \0^\top\\
 1 & 0 & 0 & \a^\top\\
 1 & 0 & 0 & \b^\top\\
 \0 & \a & \b & B
 \end {array}\right),\quad
A(\widetilde G)=\left(\begin {array} {cccc}
 0 & 0 & 0 & (\a+\b)^\top\\
 0 & 0 & 0 & \a^\top\\
 0 & 0 & 0 & \b^\top\\
 \a+\b & \a & \b & B
 \end {array}\right).
  $$

Let $s\subseteq V=V(G)\setminus\{u,v,w\}$. Then
 \begin {gather*}
\corank A(G(s))=\corank A(G(s\cup\{u\}))-1=\corank A(G(s\cup\{u,v\}))=\\
=\corank A(G(s\cup\{u,w\}))=\corank A(\widetilde G(s)),\\
\corank A(G(s\cup\{v,w\}))=\corank A(\widetilde
G(s\cup\{v,w\}))=\\
=\corank A(\widetilde G(s\cup\{u,v\}))=\corank A(\widetilde
G(s\cup\{u,w\}))=\\
=\corank A(\widetilde G(s\cup\{u,v,w\}))-1,\\
\corank A(G(s\cup\{v\}))=\corank A(\widetilde G(s\cup\{v\})),\\
\corank A(G(s\cup\{w\}))=\corank A(\widetilde
G(s\cup\{w\})),\\
\corank A(G(s\cup\{u,v,w\}))=\corank A(\widetilde G(s\cup\{u\})),
 \end {gather*}
and
   \begin {gather*}
\langle G\rangle=\sum\limits_{s\subseteq
V}\Bigl(a^{2\alpha(s)-n}(-a^2-a^{-2})^{\corank A(G(s))}+\\
+a^{2(\alpha(s)+1)-n}\bigl((-a^2-a^{-2})^{\corank
A(G(s\cup\{u\}))}+\\
+(-a^2-a^{-2})^{\corank A(G(s\cup\{v\}))}+(-a^2-a^{-2})^{\corank
A(G(s\cup\{w\}))}\bigr)+\\
+a^{2(\alpha(s)+2)-n}\bigl((-a^2-a^{-2})^{\corank
A(G(s\cup\{u,v\}))}+(-a^2-a^{-2})^{\corank A(G(s\cup\{u,w\}))}+\\
+(-a^2-a^{-2})^{\corank
A(G(s\cup\{v,w\}))}\bigr)+a^{2(\alpha(s)+3)-n}(-a^2-a^{-2})^{\corank
A(G(s\cup\{u,v,w\}))} \Bigr)=
  \end {gather*}
  \begin {gather*}
=\sum\limits_{s\subseteq
V}\Bigl(a^{2\alpha(s)-n}(-a^2-a^{-2})^{\corank
A(G(s))}\bigl(1+a^2(-a^2-a^{-2})+2a^4\bigr)+\\
+a^{2(\alpha(s)+1)-n}\bigl((-a^2-a^{-2})^{\corank
A(G(s\cup\{v\}))}+(-a^2-a^{-2})^{\corank
A(G(s\cup\{w\}))}\bigr)+\\
+a^{2(\alpha(s)+2)-n}(-a^2-a^{-2})^{\corank A(G(s\cup\{v,w\}))}+\\
+a^{2(\alpha(s)+3)-n}(-a^2-a^{-2})^{\corank
A(G(s\cup\{u,v,w\}))}\Bigr)=
 \end {gather*}
 \begin {gather*}
=\sum\limits_{s\subseteq
V}\Bigl(a^{2(\alpha(s)+2)-n}(-a^2-a^{-2})^{\corank
A(G(s))}+\\
+a^{2(\alpha(s)+1)-n}\bigl((-a^2-a^{-2})^{\corank A(G(s\cup\{v\}))}+
(-a^2-a^{-2})^{\corank A(G(s\cup\{w\}))}\bigr)+\\
+a^{2(\alpha(s)+2)-n}(-a^2-a^{-2})^{\corank A(G(s\cup\{v,w\}))}+\\
+a^{2(\alpha(s)+3)-n}(-a^2-a^{-2})^{\corank
A(G(s\cup\{u,v,w\}))}\Bigr),
  \end {gather*}
here $\alpha(s)$ is equal to the sum of the vertices labeled `$-$'
from $s$ and of the vertices labeled `$+$' from $V\setminus s$.

Analogously for $\widetilde{G}$, we get
  \begin {multline*}
\langle\widetilde{G}\rangle=\sum\limits_{s\subseteq
V}\Bigl(a^{2(\alpha(s)+2)-n}(-a^2-a^{-2})^{\corank A(\widetilde
G(s))}+\\
+a^{2(\alpha(s)+3)-n}(-a^2-a^{-2})^{\corank A(\widetilde
G(s\cup\{u\}))}+\\
+a^{2(\alpha(s)+1)-n}\bigl((-a^2-a^{-2})^{\corank A(\widetilde
G(s\cup\{v\}))}+(-a^2-a^{-2})^{\corank A(\widetilde
G(s\cup\{w\}))}\bigr)+\\
+a^{2(\alpha(s)+2)-n}\bigl((-a^2-a^{-2})^{\corank A(\widetilde
G(s\cup\{u,v\}))}+(-a^2-a^{-2})^{\corank A(\widetilde
G(s\cup\{u,w\}))}\bigr)+\\
+a^{2\alpha(s)-n}(-a^2-a^{-2})^{\corank A(\widetilde
G(s\cup\{v,w\}))}\bigr)+\\
+a^{2(\alpha(s)+1)-n}(-a^2-a^{-2})^{\corank A(\widetilde
G(s\cup\{u,v,w\}))}\Bigr)=\\
=\sum\limits_{s\subseteq
V}\Bigl(a^{2(\alpha(s)+2)-n}(-a^2-a^{-2})^{\corank A(G(s))}+\\
+a^{2(\alpha(s)+3)-n}(-a^2-a^{-2})^{\corank
A(G(s\cup\{u,v,w\}))}+\\
+a^{2(\alpha(s)+1)-n}\bigl((-a^2-a^{-2})^{\corank
A(G(s\cup\{v\}))}+(-a^2-a^{-2})^{\corank
A(G(s\cup\{w\}))}\bigr)+\\
+a^{2\alpha(s)-n}(-a^2-a^{-2})^{\corank A(\widetilde
G(s\cup\{v,w\}))}\bigl(2a^4+1+a^2(-a^2-a^{-2})\bigr)\Bigr)=\langle
G\rangle.
 \end {multline*}

   \item
Let us consider the fourth graph-move $\Omega_g 4$. Enumerate all
vertices of $G$ by numbers from $1$ to $n$ as follows. The number of
$u$ is $1$, the number of $v$ is $2$, then we enumerate the vertices
disjoint with $u$ and $v$, the vertices from $N(u)\setminus N(v)$,
the vertices from $N(v)\setminus N(u)$ and the vertices from
$N(u)\bigcap N(v)$. We define the numbers of the vertices of
$\widetilde{G}$ to be the same as the numbers of the corresponding
vertices in $G$. The adjacency matrices $A(G)$ and
$A(\widetilde{G})$ are
  \begin {gather*}
A(G)=\left(\begin {array} {cccccc}
 0 & 1 & \0^\top & \1^\top & \0^\top & \1^\top\\
 1 & 0 & \0^\top & \0^\top & \1^\top & \1^\top\\
 \0 & \0 & A_0 & A_1 & A_2 & A_3\\
 \1 & \0 & A_1 & A_4 & A_5 & A_6\\
 \0 & \1 & A_2 & A_5 & A_7 & A_8\\
 \1 & \1 & A_3 & A_6 & A_8 & A_9\\
 \end {array}\right),\\
A(\widetilde{G})=\left(\begin {array} {cccccc}
 0 & 1 & \0^\top & \1^\top & \0^\top & \1^\top\\
 1 & 0 & \0^\top & \0^\top & \1^\top & \1^\top\\
 \0 & \0 & A_0 & A_1 & A_2 & A_3\\
 \1 & \0 & A_1 & A_4 & A_5+(1) & A_6+(1)\\
 \0 & \1 & A_2 & A_5+(1) & A_7 & A_8+(1)\\
 \1 & \1 & A_3 & A_6+(1) & A_8+(1) & A_9\\
 \end {array}\right),
  \end {gather*}
where bold characters indicate row and column vectors, $(1)$ is the
matrix consisting of $1$, and $A_j$, $j=1,\dots,9$, are matrices.

Let $s\subseteq V(G)\setminus\{u,v\}$. It is not difficult to see
that
  \begin {gather*}
\corank A(G(s))=\corank A(\widetilde{G}(s\cup\{u,v\})),\\
\corank A(G(s\cup\{u,v\}))=\corank A(\widetilde{G}(s)),\\
\corank A(G(s\cup\{u\}))=\corank A(\widetilde{G}(s\cup\{u\})),\\
\corank A(G(s\cup\{v\}))=\corank A(\widetilde{G}(s\cup\{v\})),
  \end {gather*}
and $\alpha(s)$ for $G$ equals $\alpha(s\cup\{u,v\})$ for
$\widetilde{G}$, $\alpha(s\cup\{u,v\})$ for $G$ equals $\alpha(s)$
for $\widetilde{G}$, $\alpha(s\cup\{u\})$ for $G$ equals
$\alpha(s\cup\{u\})$ for $\widetilde{G}$, $\alpha(s\cup\{v\})$ for
$G$ equals $\alpha(s\cup\{v\})$ for $\widetilde{G}$.

Collecting the terms as above in the Kauffman brackets polynomial of
$G$ and $\widetilde{G}$ we immediately get $\langle
G\rangle=\langle\widetilde{G}\rangle$, which completes the proof.
  \end {enumerate}
 \end {proof}

Theorem~\ref {th:mov_reid} yields the following lemma as a
corollary.

 \begin {lem}
The Kauffman bracket polynomial of a labeled graph is invariant
under $\Omega_g 5$ graph-move.
 \end {lem}

The following theorem describes the difference of the Kauffman
bracket polynomials under $\Omega_g 1$ graph-move.

 \begin {theorem}\label {th:mov_reid_1}
The Kauffman bracket polynomial of a labeled graph is multiplied by
$(-a^{\pm3})$ under $\Omega_g 1$ graph-move. More precisely, it is
multiplied by $(-a^{-3})$ under an addition of vertex labeled `$+$'
and by $(-a^3)$ under an addition of vertex labeled `$-$'.
 \end {theorem}

 \begin {proof}
The graph $\widetilde{G}$ is obtained from $G$ by addition of one
isolated vertex $u$. Enumerate all vertices of $G$ by numbers from
$1$ to $n$, and enumerate all vertices of $\widetilde G$ by numbers
from $1$ to $n+1$ in such a way that the number of $u$ is $1$ and
the number of $v\in V(\widetilde{G})\setminus\{u\}$ is $i\geqslant2$
if and only if the number of $v$ in $G$ is $i-1$. The adjacency
matrix $A(\widetilde{G})$ is
  $$
A(\widetilde{G})=\left(\begin {array} {cc}
 0 & \0^\top\\
 \0 & A(G)
 \end {array}\right),
  $$
where bold characters indicate row and column vectors.

Let $s\subseteq V(G)$, we have
 $$
\corank A(G(s))=\corank A(\widetilde G(s))=\corank A(\widetilde
G(s\cup\{u\}))-1.
 $$

Suppose the label of $u$ is `$+$', we get
   \begin {gather*}
\langle\widetilde{G}\rangle=\sum\limits_{s\subseteq
V(G)}\Bigl(a^{2(\alpha(s)+1)-(n+1)}(-a^2-a^{-2})^{\corank
A(\widetilde G(s))}+\\
+a^{2\alpha(s)-(n+1)}(-a^2-a^{-2})^{\corank A(\widetilde
G(s\cup\{u\}))}\Bigr)=\\
=\sum\limits_{s\subseteq V(G)}a^{2\alpha(s)-n}(-a^2-a^{-2})^{\corank
A(G(s))}\bigl(a+a^{-1}(-a^2-a^{-2})\bigr)=-a^{-3}\langle G\rangle,
   \end {gather*}
here $\alpha(s)$ is equal to the sum of the vertices labeled `$-$'
from $s$ and of the vertices labeled `$+$' from $V(G)\setminus s$.

Suppose the label of $u$ is `$-$', we get
   \begin {gather*}
\langle\widetilde{G}\rangle=\sum\limits_{s\subseteq
V(G)}\Bigl(a^{2\alpha(s)-(n+1)}(-a^2-a^{-2})^{\corank A(\widetilde
G(s))}+\\
+ a^{2(\alpha(s)+1)-(n+1)}(-a^2-a^{-2})^{\corank A(\widetilde
G(s\cup\{u\}))}\Bigr)=\\
=\sum\limits_{s\subseteq V(G)}a^{2\alpha(s)-n}(-a^2-a^{-2})^{\corank
A(G(s))}\bigl(a^{-1}+a(-a^2-a^{-2})\bigr)=-a^3\langle G\rangle,
   \end {gather*}
here $\alpha(s)$ is equal to the sum of the vertices labeled `$-$'
from $s$ and of the vertices labeled `$+$' from $V(G)\setminus s$.
This completes the proof.
 \end{proof}


 \section {Graph-knots and the Jones polynomial}


In this section we define {\em graph-knots} which are analogous to
virtual knots in virtual links,  and normalise the Kauffman bracket
polynomial, in order to get the Jones polynomial.

To normalise the Kauffman bracket polynomial to an invariant (the
Jones polynomial) for the case of links (virtual links), we have to
introduce the notion of writhe number (the number of crossings of
type $\skcrro$ minus the number of crossings of type $\skcr$). For a
non-oriented knot, this number does not depend on the orientation,
whence for a link it does depend on relative orientations of the
components. Below we define graph-knots (which correspond to usual
knots --- one-component links) and define the writhe number and
hence the invariant Jones polynomial. A knot is a link having one
component, i.e.\ the set of knots is a proper subset of set of
links.

 \begin {lem}\label {lem:cor_0}
Let $G$ be a labeled graph and $\widetilde{G}$ be a graph obtained
from $G$ by some graph-move of $\Omega_g1-\Omega_g 4$. If
$\corank(A(G)+E)=0$ then $\corank(A(\widetilde G)+E)=0$, here $E$ is
the identity matrix.
 \end {lem}

 \begin {proof}
We will use the notations as the same in Theorems~\ref
{th:mov_reid},~\ref {th:mov_reid_1}.

  \begin {enumerate}
   \item
For the first graph-move $\Omega_g 1$ the assertion of the lemma
immediately follows from the view of
  $$
A(\widetilde{G})=\left(\begin {array} {cc}
 0 & \0^\top\\
 \0 & A(G)
 \end {array}\right).
  $$

   \item
Consider the second graph-move $\Omega_g 2$. Using elementary
manipulations with matrices (over $\Z_2$), we get
  \begin {gather*}
A(\widetilde{G})+E=\left(\begin {array} {ccc}
 1 & 0 & \a^\top\\
 0 & 1 & \a^\top\\
 \a & \a& A(G)+E
 \end {array}\right)\rightsquigarrow\left(\begin {array} {ccc}
 1 & 1 & \0^\top\\
 0 & 1 & \a^\top\\
 \a & \a& A(G)+E
 \end {array}\right)\rightsquigarrow\\
 \rightsquigarrow\left(\begin {array} {ccc}
 0 & 1 & \0^\top\\
 1 & 1 & \a^\top\\
 \a & \a& A(G)+E
 \end {array}\right)\rightsquigarrow\left(\begin {array} {ccc}
 0 & 1 & \0^\top\\
 1 & 0 & \0^\top\\
 \0 & \0 & A(G)+E
 \end {array}\right),
  \end {gather*}
i.e.\ $\corank(A(\widetilde{G})+E)=\corank(A(G)+E)=0$.

   \item
Consider the third graph-move $\Omega_g 3$. Using elementary
manipulations with matrices, we get
  \begin {gather*}
A(\widetilde G)+E=\left(\begin {array} {cccc}
 1 & 0 & 0 & (\a+\b)^\top\\
 0 & 1 & 0 & \a^\top\\
 0 & 0 & 1 & \b^\top\\
 \a+\b & \a & \b & B+E
 \end {array}\right)\rightsquigarrow
 \left(\begin {array} {cccc}
 1 & 1 & 1 & \0^\top\\
 0 & 1 & 0 & \a^\top\\
 0 & 0 & 1 & \b^\top\\
 \a+\b & \a & \b & B+E
 \end {array}\right)\rightsquigarrow\\
 \rightsquigarrow\left(\begin {array} {cccc}
 1 & 1 & 1 & \0^\top\\
 1 & 1 & 0 & \a^\top\\
 1 & 0 & 1 & \b^\top\\
 \0 & \a & \b & B+E
 \end {array}\right)=A(G)+E
  \end {gather*}
i.e.\ $\corank(A(\widetilde{G})+E)=\corank(A(G)+E)=0$.

   \item
For the fourth graph-move $\Omega_g 4$ we have, by using the
elementary manipulations with matrices,
  \begin {gather*}
A(\widetilde G)+E=\left(\begin {array} {cccccc}
 1 & 1 & \0^\top & \1^\top & \0^\top & \1^\top\\
 1 & 1 & \0^\top & \0^\top & \1^\top & \1^\top\\
 \0 & \0 & A_0+E & A_1 & A_2 & A_3\\
 \1 & \0 & A_1 & A_4+E & A_5+(1) & A_6+(1)\\
 \0 & \1 & A_2 & A_5+(1) & A_7+E & A_8+(1)\\
 \1 & \1 & A_3 & A_6+(1) & A_8+(1) & A_9+E\\
 \end {array}\right)\rightsquigarrow\\
\rightsquigarrow
 \left(\begin {array} {cccccc}
 1 & 1 & \0^\top & \1^\top & \0^\top & \1^\top\\
 1 & 1 & \0^\top & \0^\top & \1^\top & \1^\top\\
 \0 & \0 & A_0+E & A_1 & A_2 & A_3\\
 \1 & \0 & A_1 & A_4+E & A_5 & A_6\\
 \0 & \1 & A_2 & A_5 & A_7+E & A_8\\
 \1 & \1 & A_3 & A_6 & A_8 & A_9+E\\
 \end {array}\right),
  \end {gather*}
i.e.\ $\corank(A(\widetilde{G})+E)=\corank(A(G)+E)=0$, which
completes the proof.
  \end {enumerate}
 \end {proof}

 \begin{dfn}
A graph-link $\{G\}$ is called {\em graph-knot} if
$\corank(A(G')+E)=0$ for any representative $G'$ of the graph-link.
 \end{dfn}

Lemma~\ref {lem:cor_0} guarantees that the notion of graph-knot is
well-defined.

The following definition introduces the writhe number for a
graph-knot. That number allows us to normalise the Kauffman bracket
polynomial.

 \begin {dfn}
Set the {\em writhe number} for a labeled graph $G$ with
$\corank(A(G)+E)=0$ as follows. Enumerate all vertices of $G$ and
$B_i(G)=A(G)+E+E_{ii}$ (all elements of $E_{ii}$ except one in the
i-th column and i-th row which is equal  to $1$ are $0$) for each
vertex $v_i\in V(G)$. The writhe number $w(G)$ of $G$ is
  $$
w(G)=\sum\limits_{i=1}^n(-1)^{\corank B_i(G)}\sign v_i.
  $$
 \end {dfn}

 \begin {lem}\label {lem:wr_2_5}
The writhe number is an invariant under $\Omega_g 2 - \Omega_g 4$
graph-moves.
 \end {lem}

 \begin {proof}
Let $G$ be a graph and $w(G)$ be its writhe number. As above
$\widetilde{G}$ is the graph obtained from $G$ by some graph-move.
We have $\corank(A(\widetilde{G})+E)=0$ (Lemma~\ref {lem:cor_0}).

We will use the notation as in Theorem~\ref {th:mov_reid}.

 \begin {enumerate}
   \item
Consider the second graph-move $\Omega_g 2$. Performing the same
elementary manipulation as in Lemma~\ref {lem:cor_0} we immediately
have the equality of coranks of the vertices different from $u$ and
$v$.

As
  \begin {gather*}
w(\widetilde{G})=\sum\limits_{i=1}^{n+2}(-1)^{\corank B_i(\widetilde
G)}\sign w_i=\\
=(-1)^{\corank B_1(\widetilde G)}\sign u+(-1)^{\corank
B_2(\widetilde G)}\sign v+w(G),
  \end {gather*}
then to complete the proof in this case it is sufficient to show
that $\corank B_1(\widetilde{G})=\corank B_2(\widetilde{G})$. The
last fact is obvious.

   \item
Consider the third graph-move $\Omega_g 3$. Performing the same
elementary manipulation as in Lemma~\ref {lem:cor_0} we immediately
have the equality of coranks of the vertices different from $u$, $v$
and $w$.

As
  \begin {gather*}
w(\widetilde{G})=\sum\limits_{i=1}^n(-1)^{\corank B_i(\widetilde
G)}\sign
t_i=(-1)^{\corank B_1(\widetilde G)}\sign u+\\
+(-1)^{\corank B_2(\widetilde G)}\sign v+(-1)^{\corank
B_3(\widetilde G)}\sign v+\\
+\sum\limits_{i=4}^n(-1)^{\corank B_i(\widetilde G)}\sign
t_i=(-1)^{\corank B_1(\widetilde G)+1}+\\
+(-1)^{\corank B_2(\widetilde G)}+(-1)^{\corank B_3(\widetilde
G)}+\sum\limits_{i=4}^n(-1)^{\corank B_i(G)}\sign t_i
  \end {gather*}
and
  \begin {gather*}
w(G)=\sum\limits_{i=1}^n(-1)^{\corank B_i(G)}\sign t_i=
(-1)^{\corank B_1(G)}\sign u+\\
+(-1)^{\corank B_2(G)}\sign v+(-1)^{\corank B_3(G)}\sign v+\\
+\sum\limits_{i=4}^n(-1)^{\corank B_i(G)}\sign t_i=(-1)^{\corank
B_1(G)+1}+(-1)^{\corank B_2(G)+1}+\\
+(-1)^{\corank B_3(G)+1}+\sum\limits_{i=4}^n(-1)^{\corank
B_i(G)}\sign t_i,
  \end {gather*}
then to complete the proof in this case it is sufficient to show
that $\corank B_1(\widetilde{G})=\corank B_1(G)$, $\corank
B_2(\widetilde{G})+\corank B_3(G)=1$ and $\corank
B_3(\widetilde{G})+\corank B_2(G)=1$.

Let us prove the first and second equalities. We have:

(a)
 \begin {gather*}
B_1(\widetilde{G})=\left(\begin {array} {cccc}
 0 & 0 & 0 & (\a+\b)^\top\\
 0 & 1 & 0 & \a^\top\\
 0 & 0 & 1 & \b^\top\\
 \a+\b & \a & \b & B+E
 \end {array}\right)\rightsquigarrow
 \left(\begin {array} {cccc}
 0 & 1 & 1 & \0^\top\\
 0 & 1 & 0 & \a^\top\\
 0 & 0 & 1 & \b^\top\\
 \a+\b & \a & \b & B+E
 \end {array}\right)\rightsquigarrow\\
 \rightsquigarrow\left(\begin {array} {cccc}
 0 & 1 & 1 & \0^\top\\
 1 & 1 & 0 & \a^\top\\
 1 & 0 & 1 & \b^\top\\
 \0 & \a & \b & B+E
 \end {array}\right)=B_1(G)
  \end {gather*}
i.e.\ $\corank B_1(\widetilde{G})=\corank B_1(G)$;

(b)
 \begin {gather*}
1=\det(A(G)+E)=\det B_2(G)+\det\left(\begin {array} {cccc}
 1 & 1 & 1 & \0^\top\\
 0 & 1 & 0 & \0^\top\\
 1 & 0 & 1 & \b^\top\\
 \0 & \a & \b & B+E
 \end {array}\right),
  \end {gather*}
and
  \begin {gather*}
 \left(\begin {array} {cccc}
 1 & 1 & 1 & \0^\top\\
 0 & 1 & 0 & \0^\top\\
 1 & 0 & 1 & \b^\top\\
 \0 & \a & \b & B+E
 \end {array}\right)\rightsquigarrow
\left(\begin {array} {cccc}
 1 & 0 & 0 & \0^\top\\
 0 & 1 & 0 & \0^\top\\
 0 & 0 & 0 & \b^\top\\
 \0 & \0 & \b & B+E
 \end {array}\right),\\
B_3(\widetilde{G})=\left(\begin {array} {cccc}
 1 & 0 & 0 & (\a+\b)^\top\\
 0 & 1 & 0 & \a^\top\\
 0 & 0 & 0 & \b^\top\\
 \a+\b & \a & \b & B+E
 \end {array}\right)\rightsquigarrow
 \left(\begin {array} {cccc}
 0 & 1 & 0 & \0^\top\\
 1 & 0 & 0 & \a^\top\\
 0 & 0 & 0 & \b^\top\\
 \0 & \a & \b & B+E
 \end {array}\right)\rightsquigarrow\\
 \rightsquigarrow\left(\begin {array} {cccc}
 0 & 1 & 0 & \0^\top\\
 1 & 0 & 0 & \0^\top\\
 0 & 0 & 0 & \b^\top\\
 \0 & \0 & \b & B+E
 \end {array}\right)
  \end {gather*}
i.e.\ $1=\corank B_3(\widetilde{G})+\corank B_2(G)$.
   \item
Consider the fourth graph-move $\Omega_g 4$. Performing the same
elementary manipulation as in Lemma~\ref {lem:cor_0} we immediately
have the equality of coranks of the vertices different from $u$,
$v$.

As
  \begin {gather*}
w(\widetilde{G})=\sum\limits_{i=1}^n(-1)^{\corank B_i(\widetilde
G)}\sign w_i=(-1)^{\corank B_1(\widetilde G)}\sign u+\\
+(-1)^{\corank B_2(\widetilde G)}\sign v+
\sum\limits_{i=3}^n(-1)^{\corank B_i(\widetilde G)}\sign w_i=\\
=(-1)^{\corank B_1(\widetilde G)+1}\sign u+(-1)^{\corank
B_2(\widetilde G)+1}\sign v+\sum\limits_{i=3}^n(-1)^{\corank
B_i(G)}\sign w_i
  \end {gather*}
and
  \begin {gather*}
w(G)=\sum\limits_{i=1}^n(-1)^{\corank B_i( G)+1}\sign w_i=
(-1)^{\corank B_1(G)+1}\sign u+\\
+(-1)^{\corank B_2(G)+1}\sign v+ \sum\limits_{i=3}^n(-1)^{\corank
B_i(G)+1}\sign w_i,
  \end {gather*}
then to complete the proof in this case it is sufficient to show
that $\corank B_1(\widetilde{G})+\corank B_1(G)=1$ and $\corank
B_2(\widetilde{G})+\corank B_2(G)=1$.

Let us prove only the first equality. We have:
  \begin {gather*}
1=\det(A(G)+E)=\\
=\det B_1(G)+\det\left(\begin {array}{cccccc}
1 & 0 & \0^\top & \0^\top & \0^\top & \0^\top\\
1 & 1 & \0^\top & \0^\top & \1^\top & \1^\top\\
\0 & \0 & A_0+E & A_1 & A_2 & A_3\\
\1 & \0 & A_1 & A_4+E & A_5 & A_6\\
\0 & \1 & A_2 & A_5 & A_7+E & A_8\\
\1 & \1 & A_3 & A_6 & A_8 & A_9+E\\
 \end {array}\right),
  \end {gather*}
and
  \begin {gather*}
B_1(\widetilde{G})=\left(\begin {array} {cccccc}
 0 & 1 & \0^\top & \0^\top & \1^\top & \1^\top\\
 1 & 1 & \0^\top & \1^\top & \0^\top & \1^\top\\
 \0 & \0 & A_0+E & A_1 & A_2 & A_3\\
 \0 & \1 & A_1 & A_4+E & A_5+(1) & A_6+(1)\\
 \1 & \0 & A_2 & A_5+(1) & A_7+E & A_8+(1)\\
 \1 & \1 & A_3 & A_6+(1) & A_8+(1) & A_9+E\\
 \end {array}\right)\rightsquigarrow\\
 \rightsquigarrow
 \left(\begin {array} {cccccc}
 0 & 1 & \0^\top & \0^\top & \0^\top & \0^\top\\
 1 & 1 & \0^\top & \0^\top & \1^\top & \1^\top\\
 \0 & \0 & A_0+E & A_1 & A_2 & A_3\\
 \0 & \1 & A_1 & A_4+E & A_5 & A_6\\
 \1 & \0 & A_2 & A_5 & A_7+E & A_8\\
 \1 & \1 & A_3 & A_6 & A_8 & A_9+E\\
 \end {array}\right),
   \end {gather*}
i.e.\ $1=\corank B_1(\widetilde{G})+\corank B_1(G)$, which completes
the proof.
  \end {enumerate}
 \end {proof}

The following lemma is evident.

 \begin {lem}\label {lem:wr_1}
The writhe number is changed by $(\pm1)$ under $\Omega_g 1$
graph-move. More precisely, it is changed by $(-1)$ under an
addition of vertex labeled `$+$' and by $(+1)$ under an addition of
vertex labeled `$-$'.
 \end {lem}

 \begin {dfn}
Let $G$ be a graph-knot. Define the Jones polynomial as
$X(G)=(-a)^{-3w(G)}\langle G\rangle$.
 \end {dfn}

Theorems~\ref {th:mov_reid},~\ref {th:mov_reid_1} and Lemmas~\ref
{lem:wr_2_5},~\ref {lem:wr_1} immediately yield

  \begin{theorem}
The Jones polynomial $X(G)$ is an invariant of graph-knots.
  \end{theorem}


  \section{Minimal representatives of graph-links}


One of important problems in classical knot theory is to establish
minimal crossing number of a certain link. In late 19's century,
famous physicist and knot tabulator P.\,G.~Tait~\cite {Tait}
conjectured that alternating prime diagrams of classical links are
minimal with respect to the number of classical crossings. This
celebrated conjecture was solved only in 1987, after the notions of
the Jones polynomial and the Kauffman bracket polynomial appeared.
The first solution was obtained by Murasugi~\cite{Mur}, then it was
reproved by Thistlethwaite~\cite{Thi1}, Turaev~\cite{Turg}, and
others. Later, Thistlethwaite~\cite{Thi2} established the minimality
for a larger class of diagrams (so-called {\em adequate} diagrams).
It turns out that many results in these directions generalise for
virtual links (establishing the minimal number of classical
crossings); these results were obtained by the second named author
of the present paper. See~\cite{MBook} for the proofs and some
further generalisations and other results concerning virtual knots.

 \begin {dfn}
The difference between the leading degree and the lowest degree of
non-zero terms of the Kauffman bracket polynomial $\langle K\rangle$
is called the {\em span} of the Kauffman bracket and is denoted by
$\spn\langle K\rangle$.
 \end {dfn}

The main reason of these minimality theorems and further crossing
estimates come from the well-known Kauffman-Murasugi-Thistlethwaite
Theorem:

 \begin {theorem}
For a non-split classical link diagram $K$ on $n$ crossings we have
$\spn\langle K\rangle\leqslant 4n$, whence for alternating non-split
diagrams we have $\spn\langle K \rangle=4n$.
 \end {theorem}

Note that the span of the Kauffman bracket is invariant under all
Reidemeister moves.

In~\cite{MBook}, this theorem is generalised for virtual diagrams.
The estimate $\spn\langle K\rangle\leqslant 4n$ can be sharpened to
$\spn\langle K\rangle \leqslant 4n-4g$, where $g$ is the genus of
the atom. A nice way to reprove the latter estimate can be also
found in~\cite{DFKLS}, where the authors interpreted atoms as a
modification of Grothendieck's {\em dessins d'enfant}.

In the present section, we establish minimality of graph-link
representatives. We call a labeled graph $G$ {\em minimal} if there
is no representative of the graph-link corresponding to $G$ having
strictly smaller number of vertices than $G$ has.

 \begin {dfn}
A classical link diagram is called {\em alternating} if while
passing along every component of it we alternate undercrossings and
overcrossings.
 \end {dfn}

From the `atomic' point of view, alternating link diagrams are those
having atom genus (Turaev genus) $0$ (more precisely, diagram has
genus $0$ if it is a connected sum of several alternating diagrams).

For virtual links, we have a notion of {\em quasi-alternating}
diagram~\cite{MBook}: these are precisely diagrams obtained from
classical alternating diagrams by (detours and) virtualisations.

 \begin {dfn}
A virtual link diagram $D$ is called {\em split} if there is a
vertex $X$ of the corresponding atom $(M,\Gamma)$ such that
$\Gamma\backslash X$ is disconnected.
 \end {dfn}

Now, let us generalise the notions defined above for the case of
graph-links.

 \begin {dfn}
A labeled graph $G$ on $n$ vertices is {\em alternating} if
$k+l=n+2$, where $k$ is the number of circles in the $A$-state
$s_1$, i.e.\ $k=\corank A(G(s_1))+1$, and $l$ is the number of
circles in the $B$-state $s_2$ of $G$, i.e.\ $l=\corank
A(G(s_2))+1$.
 \end {dfn}

 \begin {dfn}
A labeled graph $G$ is {\em adequate} if the number of circles $k$
in the $A$-state is locally minimal, that is, there is no
neighbouring state for the $A$-state with $k+1$ circles, and the
same is true for the number of circles $l$ in the $B$-state.
 \end {dfn}

 \begin {rk}
This definition of the adequate diagram generalises (see,
e.g.,~\cite {Thi2}) the classical definition of the adequate
diagram: no circle of the $A$-state nor $B$-state splits into a pair
of circles after one resmoothing.
 \end {rk}

 \begin {dfn}
A labeled graph $G$ is {\em non-split} if it has no isolated
vertices.
 \end {dfn}

 \begin {dfn}
For a labeled graph $G$ let the {\em atom genus} ({\em Turaev
genus}) be $1-(k+l-n)/2$, where $k$ and $l$ are the numbers of
circles in the $A$-state and $B$-state of $G$, respectively.
 \end {dfn}

Note that this number agrees with the atom genus in the usual case:
we just use $\chi=2-2g$, where $\chi$ is the Euler characteristic,
and count $\chi$ by using the number of crossings $n$, number of
edges $2n$ and the number of $2$-cells ($A$-state circles and
$B$-state circles).

 \begin {lem}\label {lem:lmm1}
For any graph $G$ on $n$ vertices we have $\spn\langle
G\rangle\leqslant 4n-4g(G)$, here $g(G)$ is the genus of the
corresponding atom.
 \end {lem}

 \begin {proof}
Indeed, the assertion of this lemma comes from the definition of the
Kauffman bracket and the atom genus. Denote the number of circles in
the $A$-state of $G$ by $k$, and denote the number of circles in the
$B$-state of $G$ by $l$. Then the leading term of the Kauffman
bracket coming from the $A$-state has degree $n+2(k-1)$, and the
lowest term coming from the $B$-state has degree $-n-2(l-1)$. Now,
it remains to see that no other state can give a term of degree
strictly greater than that of the $A$-state. Similarly, no state
contributes a term of degree strictly smaller than that of the
$B$-state and the inequality follows.
 \end {proof}

 \begin {lem}\label {lem:lmm2}
For an adequate labeled graph $G$ on $n$ vertices we have
$\spn\langle G\rangle=4n-4g(G)$, here $g(G)$ is the genus of the
corresponding atom.
 \end {lem}

 \begin {proof}
Indeed, it is sufficient to check that the leading term coming from
the $A$-state of $G$ is not canceled by any other term coming from
another state (the argument for the lowest term coming from the
$B$-state is the same).

To do that, let us consider the term
$a^{\alpha(s)-\beta(s)}(-a^{2})^{\gamma(s)-1}$ for a state $s$. For
the $A$-state, we have $\alpha=n,\beta=0,\gamma=k$. If we switch one
crossing to the $B$-smoothing, then $\alpha$ is decreased by $1$,
$\beta$ is increased by $1$, which decreases the degree of
$a^{\alpha-\beta}$ by two. We may compensate this only if the
``number of circles'' $\gamma(s)$ (or the corresponding $\corank
A(G(s))+1$) is increased by $1$. This may happen only if there is a
state $\widetilde{s}$ adjacent to the $A$-state with $\corank
A(G(\widetilde{s}))+1=k+1$. Thus, the diagram is inadequate.
 \end {proof}

 \begin {lem}\label {lem:lmm3}
An alternating labeled graph $G$ is adequate if and only if it has
no isolated vertices.
 \end {lem}

 \begin {proof}
The direction $\Rightarrow$ is obvious. Now, assume that the diagram
$G$ is inadequate, alternating and has no isolated vertices. Denote
the number of circles of the $A$-state $s_1$ by $k$, and that of the
$B$-state by $l$. Without loss of generality, assume that there is a
state $A'$ with $(n-1)$ $A$-smoothings and one $B$-smoothing with
number of circles equal to $(k+1)$. Consider the opposite state
$B'$. Obviously, the number of circles in this state is $l-1$ (the
total number can not exceed $k+l$). Denote the vertex of $G$ where
the $A$-state differs from $A'$ by $X$. Thus, the labeled graph $G'$
obtained by changing the label of the $X$ has genus $0$, too.

Since $G$ is alternating, all single-circle states are at the same
distance from the $A$-state. On the other hand, all single-circle
states are at the same distance from $A'$. This means that these
single-circle states (as subsets of $\{1,\dots, n\}$) either both
contain $X$ or both do not contain $X$.

Assume they all contain $X$. We argue that $X$ is an isolated
vertex. Indeed, if there were a vertex $Y$ connected to $X$ then,
starting from a single-circle state containing $X$ and changing it
at $X$ and $Y$, we would get another single-circle state not
containing $X$.

This completes the proof.
 \end {proof}

Lemmas~\ref {lem:lmm1},~\ref {lem:lmm2},~\ref {lem:lmm3} together
yield the following

 \begin {theorem}\label {th:minimal}
Let $G$ be an alternating labeled graph without isolated vertex.
Then it is minimal, that is, there is not graph $G'$ with strictly
smaller number of vertices representing the same graph-link as $G$.
 \end {theorem}

 \begin {proof}
Assume the contrary. Then we have $4n=\spn\langle
G\rangle=\spn\langle G'\rangle\leqslant 4n'-4g(G')$, where $n'$ is
the number crossings of $G'$, and $g(G')$ is the atom genus (Turaev
genus) of $G'$. The inequality $n'<n$ leads to a contradiction,
which completes the proof.
 \end {proof}

Classical links are represented by $d$-diagrams; alternating links
are represented by those $d$-diagrams where all chords of one family
have sign `$+$' and all chords of the other family have sign `$-$'.
At the level of graphs, $d$-diagrams are bipartite graphs.

 \begin {dfn}
We call any bipartite graph with arbitrary labeling {\em
pseudo-classical}.
 \end {dfn}

It is easy to see that a labeled graph is alternating if and only if
it is pseudo-classical, and all labels of one subset of disjoint
vertices are `$+$' and all labels of the complement subset of
vertices (which are pairwise disjoint as well) are `$-$'.

\begin{figure} \centering\includegraphics[width=150pt]{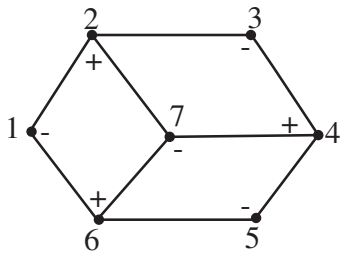}
  \caption{The graph $G_7$} \label{fig:g7}
 \end{figure}

 \begin {examp}
Consider the graph $G_{7}$ consisting of the $7$ vertices with the
following incidences. For $i,j=1,\dots, 6$, $i$ is connected to $j$
iff $i-j\equiv\pm1\pmod 6$, and $7$ is connected to $2,4,6$. Label
all even vertices by `$+$', and label all odd vertices by `$-$', see
Fig.~\ref {fig:g7}. This graph is alternating. By Theorem~\ref
{th:minimal}, $G_{7}$ is minimal. Note that this graph is not
realisable as an intersection graph of a chord diagram,
see~\cite{Bouchet}. We conjecture that the graph-link represented by
$G_{7}$ has no `realisable' representatives which correspond to
classical or virtual diagrams via chord diagrams. At least we know
that it has no such representatives with the number of crossings
less than or equal to $7$.
 \end {examp}

\begin{thebibliography}{100}

 \bibitem[BN]{BN}
D.~Bar--Natan, On the Vassiliev Knot Invariants (1995), {\em
Topology}, {\bf 34}, pp.\ 423--472.

 \bibitem[BNG]{BNG}
D.~Bar--Natan and S.~Garoufalidis, On the Melvin-Morton-Rozansky
conjecture (1996), Inv. Math. {\bf 125}, pp.\ 103--133.

 \bibitem[Bou]{Bouchet}
A.~Bouchet, Circle graph obstructions (1994), {\em J. Combinatorial
Theory B}, {\bf 60}, pp.\ 107--144.

 \bibitem[CDL]{CDL}
S.\,V.~Chmutov, S.\,V.~Duzhin, S.\,K.~Lando, Vassiliev Knot
Invariants (1994). I,II, III {\em Adv.\ Sov.\ Math.}, {\bf 21}, pp.\
117--147.

 \bibitem[DFKLS]{DFKLS}
O.~Dasbach, D.~Futer, E.~Kalfagianni, X.-S.~Lin and N.~Stoltzfus,
The Jones polynomial and graphs on surfaces (2008),  Journal of
Combinatorial Theory, B, {\bf 98} (2), pp.\ 384--399.

 \bibitem[FKM]{FKM}
R.\,A.~Fenn, L.\,H.~Kauffman, V.\,O.~Manturov, Virtual knots -
unsolved problems (2006), Fundamenta Mathematicae, N.\ 188, pp.\
293--323.

 \bibitem[Fom]{F}
A.\,T.~Fomenko, The theory of multidimensional integrable
hamiltonian systems (with arbitrary many degrees of freedom).
Molecular table of all integrable systems with two degrees of
freedom (1991), {\em Adv.\ Sov.\ Math.}, {\em 6}, pp.\ 1--35.

 \bibitem[GPV]{GPV}
M.~Goussarov, M.~Polyak, and O.~Viro // Topology. 2000. V.\ 39, pp.\
1045--1068.

 \bibitem[Jon]{Jones}
V.\,F.\,R.~Jones, A polynomial invariant for links via Neumann
algebras (1985), {\em Bull.\ Amer.\  Math.\  Soc.}, {\bf 129}, pp.\
103--112.

 \bibitem[Ka1]{KaV}
L.\,H.~Kauffman, Virtual knot theory, Eur. J. Combinatorics. 1999.
V.\ 20, N.\ 7, pp.\ 662--690.

 \bibitem[Ka2]{KauffmanBracket}
L.\,H.~Kauffman, State Models and the Jones Polynomial (1987), {\em
Topology}, {\bf 26}, pp.\ 395--407.

 \bibitem[Kh]{Kho}
M.~Khovanov, A categorification of the Jones polynomial (1997), {\em
Duke Math. J},{\bf 101} (3), pp.\ 359--426.

 \bibitem[KhR1]{KhR1}
M.~Khovanov, L.~Rozansky, Matrix Factorizations and Link Homology,
arxiv.math: GT/0401268.

 \bibitem[KhR2]{KhR2}
M.~Khovanov, L.~Rozansky, Matrix Factorizations and Link Homology
II, arxiv.math: GT/050506.

 \bibitem[Kup]{Kup}
G.~Kuperberg, What is a Virtual Link? (2002), www.arxiv.org,
math-GT$\slash$ 0208039, {\em Algebraic and Geometric Topology},
2003, {\bf 3}, 5 87-591.

 \bibitem[Lando]{Lando}
S.\,K.~Lando, $J$-invariants of ornaments and framed chord diagrams
(2006), {\em Functional Analysis and Its Applications}, {\bf 40}
(1), pp.\ 1--13.

 \bibitem[Ma1]{Izv}
V.\,O.~Manturov, Khovanov homology of virtual knots with arbitrary
coefficients), Russ.\ Acad.\ Sci.\ Math.\ Izvestiya, 2007, Vol.\ 71,
N.\ 5, pp.\ 967-999.

 \bibitem[Ma2]{Arxiv}
V.\,O.~Manturov, Embeddings of four-valent framed graphs into
2-surfaces, arxiv.math: GT/0804.4245

 \bibitem[Ma3]{AtomsandKnots}
V.\,O.~Manturov, Bifurcations, Atoms, and Knots (2000), Moscow
Univ.\ Math.\ Bull. {\bf 1}, pp.\ 3--8.

 \bibitem[Ma4]{MBook}
V.\,O.~Manturov, {\em Teoriya Uzlov} (Knot Theory), (Moscow-Izhevsk,
RCD), 2005 (512 pp.).

 \bibitem[Mur]{Mur}
K.~Murasugi, The Jones polynomial and classical conjectures in knot
theory (1987), {\em Topology} {\bf 26}, pp.\ 187-194.

 \bibitem[Oht]{Oht}
T.~Ohtsuki, Quantum Invariants (2002). A Study of Knots,
3-Manifolda, and their Sets, {\em Singapore: World Scientific}.

 \bibitem[Sob]{Soboleva}
E.~Soboleva, Vassiliev Knot Invariants Coming from Lie Algebras and
$4$-Invariants (2001),  {\em Journal of Knot Theory and Its
Ramifications}, {\bf 10} (1), pp.\ 161--169.

 \bibitem[Tait]{Tait}
P.\,G.~Tait, On knots (1898), in {\em Scientific paper I}, (London:
Campbridge University Press), pp.\ 273--317.

 \bibitem[Thi1]{Thi1}
M.\,B.~Thistlethwaite, Kauffman polynomial and alternating links
(1988), {\em Topology}, {\bf 27}, pp.\ 311-–318

 \bibitem[Thi2]{Thi2}
M.\,B.~Thistlethwaite, On the Kauffman polynomial of an adequate
link (1988), Invent.\ Math. {\bf 93} (2), pp.\ 285--296.

 \bibitem[TZ]{TZ}
L.~Traldi, L.~Zulli, A bracket polynomial for graphs, arxiv.math:
GT/0808.3392.

 \bibitem[Tu]{Turg}
V.~Turaev, A simple proof of the Murasugi and Kauffman theorems on
alternating links (1987), {\em L'Enseignement Math\'ematique}, {\bf
33}, pp.\ 203--225.

 \bibitem[Weh]{Weh}
S.~Wehrli, Khovanov homology and Conway mutations (2003),
arxiv.math: GT/0301312.
  \end{thebibliography}
 \end {document}